\documentclass[11pt]{article}

\usepackage{amssymb,latexsym}
\usepackage{amsmath}
\usepackage{citesort}
\usepackage{eucal}

\title{$^*$-Ideals and Formal Morita Equivalence of $^*$-Algebras}

\author{{\bf     
         Henrique Bursztyn\thanks{henrique@math.berkeley.edu}
         \thanks{Research supported by a fellowship from CNPq, Grant
           200481/96-7.} 
        } \\[0.5cm]
        Department of Mathematics\\
        UC Berkeley\\
        94720 Berkeley, CA, USA
        \\[1cm]
        {\bf
         Stefan
         Waldmann\thanks{Stefan.Waldmann@ulb.ac.be}
         \thanks{Research 
         supported by the Communaut\'e fran\c caise de Belgique,
         through an Action de Recherche Concert\'ee de la Direction de la
         Recherche Scientifique.}
        } \\[0.5cm]
        D{\'{e}}partement de Math{\'{e}}matique \\
        Campus Plaine, C. P. 218 \\
        Boulevard du Triomphe \\
        B-1050 Bruxelles \\
        Belgique 
       }

\date{May 2000}

\newcommand{\axiom}[1] {\textbf{{#1}}}

\newcommand{\im} {{\mathrm i}}

\newcommand{\cc} [1]     {\overline {{#1}}}
\newcommand{\supp}       {\mathop{{\mathrm {supp}}}}

\newcommand{\id}         {{\mathsf {id}}}

\newcommand{\End}        {{\mathsf {End}}}
\newcommand{\ring}[1]    {{\mathsf {{#1}}}}

\newcommand{\SP} [1]     {{\left\langle {{#1}} \right\rangle}}
\newcommand{\Bounded}    {{\mathfrak B}}
\newcommand{\Unit}       {\mathsf {1}}
\newcommand{\cl}         {\mathrm{cl}}
\newcommand{\snd}        {\mathrm{snd}}

\newcommand{\Jmin} {\mathcal J_{\mathrm {min}}}

\newcommand{\BXA} {{\sideset{_{\scriptscriptstyle \mathcal B}}
                            {_{\scriptscriptstyle\mathcal A}}
                            {\operatorname{\mathfrak X}}}}
\newcommand{\AXB} {{\sideset{_{\scriptscriptstyle \mathcal A}}
                            {_{\scriptscriptstyle\mathcal B}}
                            {\operatorname{\overline {\mathfrak X}}}}}

\newcommand{\tildeK} {\widetilde{\mathfrak K}}

\newcommand{\piA} {\pi_{\scriptscriptstyle \mathcal A}}
\newcommand{\piB} {\pi_{\scriptscriptstyle \mathcal B}}

\newcommand{\SPA}[1] {\left\langle{{#1}}
                       \right\rangle_{\scriptscriptstyle \mathcal A}}
\newcommand{\SPB}[1] {\left\langle{{#1}}
                       \right\rangle_{\scriptscriptstyle \mathcal B}}

\newcommand{\BSP}[1] {{}_{\scriptscriptstyle \mathcal B}\!
                       \left\langle {{#1}} \right\rangle}

\newcommand{\SPH}[1] {\left\langle{{#1}}
                       \right\rangle_{\scriptscriptstyle \mathfrak H}}

\newcommand{\SPKT}[1] {\left\langle{{#1}}
                       \right\rangle_{\scriptscriptstyle \widetilde{\mathfrak K}}}

\newcommand{\SPK}[1] {\left\langle{{#1}}
                       \right\rangle_{\scriptscriptstyle \mathfrak K}}

\newcommand{\srepA}  {{{}^*\textrm{-}\mathsf{rep}(\mathcal A)}}

\newcommand{\sRepA}  {{{}^*\textrm{-}\mathsf{Rep}(\mathcal A)}}

\newcommand{\srepAJ} {{{}^*\textrm{-}\mathsf{rep}(\mathcal A\big/\mathcal J)}}
\newcommand{\sRepAJ} {{{}^*\textrm{-}\mathsf{Rep}(\mathcal A\big/\mathcal J)}}

\newcommand{\srepDSA} {{{}^*\textrm{-}\mathsf{repDS}(\mathcal A)}}
\newcommand{\sRepDSA} {{{}^*\textrm{-}\mathsf{RepDS}(\mathcal A)}}

\newcommand{\RieffelX} {\mathfrak R_{\scriptscriptstyle \mathfrak X}}
\newcommand{\RieffelXc} {\mathfrak R_{\scriptscriptstyle \overline{\mathfrak X}}}

\newcommand{\PhiX}{\Phi_{\mathfrak X}}
\newcommand{\PhiccX}{\Phi_{\cc{\mathfrak X}}}

\newtheorem{lemma} {Lemma} [section]
\newtheorem{proposition} [lemma] {Proposition}
\newtheorem{theorem} [lemma] {Theorem}
\newtheorem{corollary} [lemma] {Corollary}
\newtheorem{definition}[lemma] {Definition}

\newtheorem{remark}[lemma]{Remark}

\newenvironment{proof}{\small{\sc Proof:}}{{\hspace*{\fill} $\square$\\}}

\numberwithin{equation}{section}

\begin{document}

\maketitle

\begin{abstract}
Motivated by deformation quantization, we introduced in an earlier work
the notion of formal Morita equivalence in the category of $^*$-algebras
over a ring $\ring C$ which is the quadratic extension by $\im$ of an
ordered ring $\ring R$. The goal of the present paper is twofold.
First, we clarify the relationship between formal Morita equivalence,
Ara's notion of Morita $^*$-equivalence of rings with involution,
and strong Morita equivalence of $C^*$-algebras. Second, in the general setting 
of $^*$-algebras over $\ring C$,
we define  `closed' $^*$-ideals as the ones occuring as kernels of 
$^*$-representations of these algebras on pre-Hilbert spaces. These ideals
form a lattice which we show is invariant under formal Morita equivalence.
This result, when applied to Pedersen ideals of $C^*$-algebras, recovers the
so-called Rieffel correspondence theorem. The  triviality of the
minimal element in the lattice of closed ideals, called the `minimal ideal', 
is also a formal Morita invariant and this fact can be used to describe a large
class of examples of $^*$-algebras over $\ring C$ with equivalent representation
theory but which are not formally Morita equivalent. We finally compute the
closed $^*$-ideals of some $^*$-algebras arising in differential geometry.

\end{abstract}

\newpage

\tableofcontents


%
%

\section{Introduction}
\label{IntroSec}

The concept of formal deformation quantization, first introduced in
\cite{BFFLS78}, has been used successfully to
construct and classify the quantum observable algebras of physical
systems in terms of classical data, generally given by 
a symplectic or Poisson manifold (see
\cite{DL83b,Fed94a,OMY91,Kon97b,NT95a,NT95b,BCG97,WX98} for existence
and classifications and \cite{Wei94,Ste98,Gutt2000} for recent
reviews). One is naturally led to investigate representations of these
algebras, as different representations can be used to describe
different physical situations. It is then important to decide what kind of
representations (or category of modules) of the observable algebras 
have  physical meaning and should be considered. It turns out that the
usual $^*$-representation theory of $C^*$-algebras provides a
guide for the formulation of representations of quantum observable algebras
arising in deformation quantization. In particular,
one can always consider star-products compatible with the natural
$^*$-involution on functions (given by complex conjugation), and we note that
the underlying ring of real formal power series $\mathbb R[[\hbar]]$ has
a natural order structure. Combining these two facts, one can define
positive $\mathbb C[[\hbar]]$-linear functionals on the star-product algebras 
 and $^*$-representations on pre-Hilbert spaces over $\mathbb C[[\hbar]]$
analogous to the $C^*$-algebra case. For example, 
every positive linear functional gives rise to a GNS construction of a 
$^*$-representation (see \cite{BW98a} for a 
detailed exposition of these ideas). 
These notions have yielded physically interesting representations 
and, in fact, most known operator representations are of this type, see
e.g.~\cite{BNW98a,BNW99a,BNPW98,Wal2000a}.

With this motivation, we started in \cite{BuWa99a,BuWa99b} a more systematic
study of $^*$-algebras over $\ring C = \ring R(\im)$, where $\ring R$
is an arbitrary (commutative, unital) ordered ring and $\ring C$ 
is its quadratic extension
by $\im$, with $\im^2 = -1$, and their representations on pre-Hilbert
spaces over $\ring C$. We note that not only star-product algebras
over $\mathbb C[[\hbar]]$ but also arbitrary $^*$-algebras over $\mathbb C$
(such as $^*$-algebras of unbounded operators, see
e.g.~\cite{Schmue90}) fit into this framework.
In this setting, one can consider a purely algebraic
version of Rieffel induction and define the notion of formal Morita
equivalence \cite{BuWa99a}, analogous to $C^*$-algebras 
\cite{Rief74b,Rief74,Land98}. In the present paper, we continue 
the investigation of $^*$-representations of $^*$-algebras and formal
Morita equivalence.

First, we clarify the relationship between formal Morita equivalence,
Ara's notion of Morita $^*$-equivalence of rings with involution
\cite{Ara99a,Ara99b} and strong Morita equivalence of $C^*$-algebras
\cite{Rief82}. For non-degenerate and idempotent
$^*$-algebras over $\ring C$, we show that formal Morita equivalence implies
Morita $^*$-equivalence. In the case where the $^*$-algebras 
are given by Pedersen ideals of $C^*$-algebras \cite{Peder79}, 
we show that, in fact,
formal Morita equivalence and Morita $^*$-equivalence are equivalent. It 
then follows immediately from Ara's recent algebraic characterization
of strong Morita equivalence of $C^*$-algebras in terms of their Pedersen
ideals \cite[Thm.~2.4]{Ara99b} that two 
$C^*$-algebras are strongly Morita equivalent if and only if
their Pedersen ideals are formally Morita equivalent.

Second, we discuss $^*$-ideals in $^*$-algebras over $\ring C$, their
behavior under formal Morita equivalence and how they can be used to 
describe particular properties of star-product algebras which 
general $^*$-algebras over $\ring C$ may lack. In a $^*$-algebra over
$\ring C$, we define an ideal to be `closed' if it occurs as the
kernel of some $^*$-representation  on a pre -Hilbert space over
$\ring C$. These ideals form a lattice, which we show is invariant under
formal Morita equivalence. In the case of $C^*$-algebras, an ideal is closed
in this algebraic sense if and only if it is norm closed and the aforementioned
invariance reduces to the so-called Rieffel correspondence theorem
(see \cite[Prop.~3.24]{RaeWil98}). We note that $^*$-algebras having a ``large''
number of positive linear functionals possess nice features which
do not hold in general. For instance, they admit faithful representations
on pre-Hilbert spaces (see \cite[Prop.~2.8]{BuWa99a}). We observe that not only
$C^*$-algebras but also star-product algebras arising in deformation
quantization have this property \cite[Prop.~5.3]{BuWa99b}. Admitting a
faithful representation is equivalent to the triviality of the minimal
closed ideal (defined as the intersection of all closed $^*$-ideals),
which we also show is a formal Morita invariant.  We use this invariant
to describe a large class of examples of $^*$-algebras
having equivalent categories of (strongly non-degenerate)
$^*$-representations on pre-Hilbert spaces but which are not formally
Morita equivalent, generalizing \cite[Cor.~5.20]{BuWa99a}. (It was
shown in \cite[Thm.~5.10]{BuWa99a} that formal Morita equivalence
implies the equivalence of categories of strongly non-degenerate
$^*$-representations.)

Finally, we apply these general algebraic notions to the particular case of
the $^*$-algebra $C^\infty(M)$ of complex-valued smooth functions on a manifold
$M$. We show that closed $^*$-ideals in this example are given by the ideals of 
functions vanishing on closed subsets of $M$, and hence agree with the closed
ideals with respect the (weak) $C^0$-topology (\cite[Thm.~3]{Whit48}).
We also show that $C^\infty(M)$
is formally Morita equivalent to the $^*$-algebra of smooth sections of the
endomorphism bundle of any Hermitian vector bundle $E$ over $M$. These
facts are not surprising since their analogues in the continuous 
($C^*$-algebraic) category (in the sense of strong Morita equivalence)
are well-known. They are, however, conceptually important since they show that
$^*$-algebras in the smooth category can be treated independently through
general algebraic methods (we note that $M$ is \emph{not} assumed to be compact
and hence the $^*$-algebras in questions are \emph{not} pre-$C^*$-algebras).
We point out that this discussion is the `classical' starting point in
the investigation of these ideas in the framework of formal deformation 
quantization, see also \cite{LR88} for the case of a trivial bundle.

The paper is organized as follows. In Section~\ref{PrelimSec}, we
recall the notions of ordered rings, $^*$-algebras, and their
$^*$-representations, including the GNS construction (details can be found
in \cite{BW98a,BuWa99a}). 
In Section \ref{C*worldSec}, we describe how the notions of formal
Morita equivalence, Morita $^*$-equivalence and strong Morita equivalence are
related.
We define closed $^*$-ideals
in Section~\ref{MinimalSec} and discuss basic properties of the
minimal closed $^*$-ideal and its relation to GNS representations and
positive functionals. In Section~\ref{MoritaSec}, we focus on
algebraic Rieffel induction and formal Morita equivalence of
$^*$-algebras and prove that the lattice of closed $^*$-ideals is a
formal Morita invariant. Section~\ref{FunctionSec} contains the
description of closed $^*$-ideals of $C^\infty(M)$.
As a $^*$-algebra with non-trivial minimal ideal we consider the
sections of the Grassmann bundle over $M$. We also show in this 
section that for any Hermitian vector bundle $E$ over a manifold $M$,
$C^\infty(M)$ and $\Gamma^\infty(\End(E))$ are formally Morita
equivalent, as well as $C^\infty_0 (M)$ and $\Gamma^\infty_0 (\End(E))$.

%
%

\section{Preliminaries}
\label{PrelimSec}

In this section, we just recall some basic facts on ordered rings,
pre-Hilbert spaces and $^*$-algebras to establish our notation
(see \cite{BW98a,BuWa99a} for details).

Let $\ring R$ be an associative  commutative ring with $1 \ne 0$ which
is \emph{ordered}, i.e. equipped with a distinguished subset 
$\ring P$, the positive elements, such that $\ring R$ is the disjoint
union $-\ring P \cup \{0 \} \cup \ring P$ and one has
$\ring P + \ring P, \ring P \cdot \ring P \subseteq \ring P$. We write
$a > b$ if $a-b \in \ring P$ etc. In particular $1 > 0$, and $\ring R$
has characteristic zero and no zero divisors. Next we consider the
canonical quadratic ring extension $\ring C = \ring R(\im)$ where
$\im^2 = -1$. Then $\ring C$ becomes an associative commutative ring
with $1 \ne 0$ which has again characteristic zero and no zero
divisors. Complex conjugation $z \mapsto \cc z$ is defined as usual
for $z \in \ring C$ and the embedding 
$\ring R \hookrightarrow \ring C$ determines the real elements of
$\ring C$, i.e. those with $z = \cc z$. Finally, note that 
$\cc z z \ge 0$ and $\cc z z = 0$ if and only if $z = 0$. Examples of
ordered rings are of course $\mathbb Z$, $\mathbb R$, and the formal
power series $\mathbb R[[\lambda]]$, which is the important example for
formal deformation quantization.

A \emph{pre-Hilbert space} over $\ring C$ is a $\ring C$-module
$\mathfrak H$ with positive-definite Hermitian product 
$\SP{\cdot,\cdot}: \mathfrak H \times \mathfrak H \to \ring C$ which
is linear in the second argument (physicists' convention). We
immediately have the Cauchy-Schwarz inequality
$\SP{\phi,\psi}\SP{\psi,\phi} \le \SP{\phi,\phi} \SP{\psi,\psi}$
for all $\phi, \psi \in \mathfrak H$. For a pre-Hilbert space
$\mathfrak H$, we say that $A \in \End_{\ring C} (\mathfrak H)$ has a
(necessarily unique) adjoint $A^* \in \End_{\ring C} (\mathfrak H)$ if
for all $\phi, \psi \in \mathfrak H$ the equality 
$\SP{A^*\phi, \psi} = \SP{\phi, A\psi}$ holds. 
One defines 
$\Bounded (\mathfrak H) := \{ A \in \End_{\ring C} \; | \; A^*$ exists
$\}$ and easily finds that $\Bounded (\mathfrak H)$ is an associative
algebra over $\ring C$ with unit element $\id_{\mathfrak H}$ and a
$^*$-involution $A \mapsto A^*$. More generally, one defines a
\emph{$^*$-algebra} $\mathcal A$ to be an associative algebra over
$\ring C$ equipped with a $^*$-involution, i.e. an antilinear
involutive anti-automorphism. Hermitian, normal, isometric and unitary
elements in $\mathcal A$ are defined as usual. An 
\emph{approximate identity} consists of a filtration 
$\mathcal A = \bigcup_{\alpha \in I} \mathcal A_\alpha$ by 
$\ring C$-submodules $\{\mathcal A_\alpha\}_{\alpha \in I}$ indexed by
a directed set $I$ and elements $\{E_\alpha\}_{\alpha \in I}$ such
that $E_\alpha^* = E_\alpha$, 
$E_\alpha E_\beta = E_\alpha = E_\beta E_\alpha$ for $\alpha < \beta$
and $E_\alpha A = A = A E_\alpha$ for all $A \in \mathcal A_\alpha$
and $\alpha \in I$. 
A linear functional $\omega: \mathcal A \to \ring C$ is called
\emph{positive} if $\omega(A^*A) \ge 0$ for all 
$A \in \mathcal A$. Then one has the Cauchy-Schwarz inequality 
$\omega (A^*B) \cc{\omega (A^*B)} \le \omega(A^*A) \omega(B^*B)$ as
well as the `almost reality' property $\omega(A^*B) = \cc{\omega(B^*A)}$, which
implies reality $\omega(A^*) = \cc{\omega(A)}$
if e.g. $\mathcal A$ has an approximate identity. 
Using positive functionals, we are able to define positive algebra
elements as well: by $\mathcal A^{++}$ we denote the
\emph{algebraically positive elements} of the form 
$A = a_1 A_1^*A_1 + \cdots + a_n A_n^*A_n$ where $a_i > 0$ and 
$A_i \in \mathcal A$. By $\mathcal A^+$ we denote the \emph{positive}
elements of $\mathcal A$, where $A \in \mathcal A^+$ if $A$ is Hermitian and
$\omega(A) \ge 0$ for all positive linear functionals 
$\omega: \mathcal A \to \ring C$. We note that if $\mathcal A$ is a 
$C^*$-algebra, then these notions of positivity agree with the standard ones.
Moreover, in this case $\mathcal A^+ = \mathcal A^{++}$. These notions also
yield the expected results for $C^\infty(M)$, see \cite[App.~B]{BuWa99a}.

A \emph{$^*$-representation} of
a $^*$-algebra $\mathcal A$ on a pre-Hilbert space $\mathfrak H$ 
is a $^*$-homomorphism 
$\pi: \mathcal A \to \Bounded (\mathfrak H)$. Such a
$^*$-representation is called \emph{non-degenerate} if 
$\pi(A)\phi = 0$ for all $A \in \mathcal A$ implies $\phi = 0$ and
\emph{strongly non-degenerate} if the $\ring C$-span of all
$\pi(A)\phi$ with $A \in \mathcal A$ and $\phi \in \mathfrak H$
coincides with $\mathfrak H$. In general, strong non-degeneracy
implies non-degeneracy and if $\mathcal A$ has a unit both notions
agree and are equivalent to the condition 
$\pi(\Unit) = \id_{\mathfrak H}$.

If $\omega: \mathcal A \to \ring C$ is a positive linear functional,
then 
$\mathcal J_\omega := \{A \in \mathcal A \;|\; \omega(A^*A) = 0\}$
is a left ideal in $\mathcal A$, the so-called Gel'fand ideal and
hence $\mathfrak H_\omega := \mathcal A \big/ \mathcal J_\omega$ is a
$\mathcal A$-left module. Moreover, $\mathfrak H_\omega$ becomes a
pre-Hilbert space via $\SP{\psi_B,\psi_C} := \omega(B^*C)$ such that
the left action $\pi_\omega(A)\psi_B := \psi_{AB}$ is a
$^*$-representation, the so-called GNS representation. Here 
$\psi_B \in \mathfrak H_\omega$ denotes the equivalence class of 
$B \in \mathcal A$. This GNS construction will play a crucial role
in this paper (see \cite{BW98a} for an extensive treatment).

Let us  recall the notion of algebraic Rieffel
induction, as discussed in \cite{BuWa99a}. Let $\mathcal A$ be a $^*$-algebra
over $\ring C$ and let $\mathfrak{X}$ be a right 
$\mathcal A$-module ({\it all} modules over $^*$-algebras considered here
will be assumed to have a compatible $\ring C$-module structure).
An \emph{$\mathcal A$-valued inner product} is a map $\SPA{\cdot,\cdot}:
\mathfrak{X}\times \mathfrak{X} \longrightarrow \mathcal A$, which is 
$\ring C$-linear and $\mathcal A$-right linear in the second argument
and satisfies $\SPA{x,y} = \SPA{y,x}^*$, for all $x,y \in \mathfrak{X}$.
We call this inner product \emph{positive semi-definite} if
$\SPA{x,x} \in \mathcal{A}^+$ and \emph{positive definite} if in addition
$\SPA{x,x}=0$ implies $x=0$. We say it is \emph{full} if  the 
$\ring C$-span of all $\SPA{x, y}$ is $\mathcal A$.
A (full, positive semi-definite) $\mathcal A$-valued inner product
on a left $\mathcal A$-module is defined analogously, but with linearity
properties (with respect to $\ring C$ and $\mathcal A$) in the first argument.
If $\mathcal B$ is another $^*$-algebra over $\ring C$, and $\BXA$ is a 
$(\mathcal B$-$\mathcal A)$-bimodule with an $\mathcal A$-valued inner product,
we call the $\mathcal B$-action  \emph{compatible} if $\SPA{B\cdot x,y}=
\SPA{x,B^* \cdot y}$ for all $B \in \mathcal{B}$ and 
$x,y \in \BXA$.

Let  $\BXA$ be a $(\mathcal B$-$\mathcal A)$-bimodule with positive
semi-definite $\mathcal A$-valued inner product and compatible 
$\mathcal B$-action. If 
$\piA$ is a $^*$-representation of $\mathcal A$ on $\mathfrak H$, then
one can consider the $\ring C$-module 
$\tildeK := \BXA \otimes_{\mathcal A} \mathfrak H$, where the 
$\mathcal A$-balanced tensor product is defined using $\piA$. Then
$\tildeK$ is a $\mathcal B$-left module and by 
$\SPKT{x\otimes\psi, y\otimes\phi} 
:= \SPH{\psi, \piA(\SPA{x,y})\phi}$
one obtains a Hermitian inner product on $\tildeK$. We say that the
bimodule $\BXA$ satisfies the property \textbf{P} if this Hermitian
product is positive semi-definite for all $^*$-representations 
$(\piA, \mathfrak H)$ of $\mathcal A$, see 
\cite[Sect.~3 and 4]{BuWa99a} for various sufficient conditions
guaranteeing \textbf{P}. In this case one can quotient
out the length zero vectors $\tildeK^\bot$ to obtain a pre-Hilbert
space $\mathfrak K = \tildeK \big/ \tildeK^\bot$. It is easy to check
that the left action of $\mathcal B$ passes to the quotient and yields
a $^*$-representation $\piB$ on $\mathfrak K$. The whole procedure,
called \emph{algebraic Rieffel induction}, is functorial and the functor is
denoted by $\RieffelX$. Given a 
$(\mathcal B$-$\mathcal A)$-bimodule with a positive semi-definite
$\mathcal B$-valued inner product and compatible $\mathcal A$-action, we
can consider the corresponding conjugated 
$(\mathcal A$-$\mathcal B)$-bimodule, denoted by 
$\overline{\mathfrak{X}}$, which is just $\mathfrak{X}$ as additive group,
but with $\mathcal A$ and $\mathcal B$ acting by adjoints and $\ring C$ by
complex conjugates. We then say that 
$\BXA$ satisfies property \axiom{Q} if $\overline{\mathfrak{X}}$ satisfies
\axiom{P} as above (see \cite[Sect.~5]{BuWa99a}).

Let us now consider formal Morita equivalence, as in \cite{BuWa99a}. 
Suppose $\BXA$ is a $(\mathcal B$-$\mathcal A)$-bimodule, equipped with 
full positive semi-definite $\mathcal A$-valued and $\mathcal B$-valued
inner products, compatible $\mathcal A$ and $\mathcal B$ actions and so that
properties \axiom{P} and \axiom{Q} are satisfied. 
We also assume that the inner products are \emph{compatible}, in the sense that
$x\cdot \SPA{y,z} = \BSP{x,y}\cdot z$ for all $x,y,z \in \mathfrak{X}$.
We call this object an
\emph{equivalence bimodule} and say that two $^*$-algebras $\mathcal A$ and
$\mathcal B$ are \emph{formally Morita equivalent} if there exists a
$(\mathcal B$-$\mathcal A)$-equivalence bimodule. We recall that 
if the actions of 
$\mathcal A$ and $\mathcal B$ on $\BXA$ are in addition strongly
non-degenerate, then $\BXA$ is called a \emph{non-degenerate equivalence
bimodule} and one can show that in this case 
that $\mathcal A$ and $\mathcal B$ have equivalent categories
of (strongly non-degenerate) $^*$-representations, with the functors
$\RieffelX$ and $\RieffelXc$ defining this equivalence
(\cite[Thm.~5.10]{BuWa99a}). As for notation, we will write
$\srepA$ for the category of $^*$-representations
of $\mathcal A$ (usually considered with isometric interwiners as morphisms)
and $\sRepA$ for the category of strongly non-denerate $^*$-representations.

A $^*$-algebra over $\ring C$ has 
\emph{sufficiently many positive linear functionals} if for every
non-zero Hermitian element $H$ one finds a 
positive linear functional $\omega$ with $\omega (H) \ne 0$. We recall the 
following important fact: a $^*$-algebra $\mathcal
A$ \emph{with approximate identity} has sufficiently many positive linear
functionals if and only if $\mathcal A$ has a faithful
$^*$-representation. As a consequence, these algebras present some nice
properties, resembling $C^*$-algebras. For instance, 
they are torsion-free ($zA=0$ implies $A=0$ or $z=0$, for 
$z \in \ring C$), $A^*A = 0$ implies $A = 0$ 
and $H^n = 0$ for a normal element $H$ implies $H = 0$ 
(\cite[Prop.~2.8]{BuWa99a}). We observe that the assumption on the existence
of an approximate identity can be weakened in particular cases. Indeed, if 
$\mathcal A$ is a $C^*$-algebra or a star-product algebra, 
then it always has sufficiently
many positive linear functionals but it does not admit an algebraic 
approximate identity in general. Nevertheless,
these algebras always have faithful representations and all the aforementioned
nice properties. We observe that in both cases, however, a topological
approximate identity exists (with respect to the norm topology in
the case of $C^*$-algebras and with respect to the
$\lambda$-adic topology for star-product algebras, see 
\cite[Sect.~8]{BuWa99a}).

On the other hand, the existence of sufficiently many positive linear 
functionals, without any further requirement, is generally not enough
to guarantee the existence of faithful representations and desirable
algebraic properties. For instance, consider $\ring C = \ring R(\im)$
with its natural $\ring C$-module structure. Define the $^*$-algebra 
$\mathcal A$ to be this $\ring C$-module, with the natural $^*$-involution
but with zero multiplication. Then \emph{all} linear functionals on
$\mathcal A$ are positive and $\mathcal A$ has sufficiently many positive linear
functionals. Note that $\mathcal A$ does not admit faithful representations
(since we have $A\neq 0$ with $A^*A=0$). For an example of a $^*$-algebra
\emph{without} sufficiently many positive linear functionals, let
$\bigwedge^\bullet (\ring C^n)$ be the Grassmann algebra
over $\ring C^n$, equipped with a $^*$-involution given by $e_i^* := e_i$, 
where $e_1, \ldots, e_n$ is the canonical basis of $\ring C^n$.
Denote by $\bigwedge^+ (\ring C^n)$ the ideal of
elements with positive degree and write $\bigwedge^\bullet (\ring C^n) 
= \ring C\Unit \oplus \bigwedge^+ (\ring C^n)$. Then every
positive linear functional  $\omega: \bigwedge^\bullet (\ring C^n) 
\to \ring C$ has to satisfy $\omega (\Unit) \ge 0$ and 
$\omega|_{\bigwedge^+ (\ring C^n)} = 0$. Hence, up to normalization, 
there is only one positive functional. In particular, every
(Hermitian) element in $\bigwedge^+ (\ring C^n)$ is positive but 
$\bigwedge^\bullet (\ring C^n)$ 
does not have sufficiently many positive linear functionals for 
$n \ge 1$. Note that this $^*$-algebra is nevertheless
torsion-free. For an example 
with torsion, set $\ring R = \mathbb Z$ and $\ring C = \mathbb Z \oplus 
\im \mathbb Z$ and consider $\mathbb Z_2 = \{0, \Unit\}$. 
Then $\mathbb Z_2$ becomes a unital $^*$-algebra over $\ring C$ by 
$\im \cdot \Unit = \Unit$ and $\Unit^* = \Unit$. 
Note that $\mathbb Z_2$ is not torsion-free, since for instance 
$2 \cdot \Unit = 0$, and there are no 
$\ring C$-linear functionals $\omega: \mathbb Z_2 \to \ring C$
except for the zero functional.

Thus the general
category of $^*$-algebras over $\ring C$ seems too wide for the
interpretation of these algebras as `observable algebras'. We will
discuss in Sections \ref{MinimalSec} and \ref{MoritaSec} how the ideal
strucuture of these algebras  help characterizing the more `observable
algebra'-like objects.

%
%

\section{Formal Morita equivalence and $C^*$-algebras}
\label{C*worldSec}

The aim of this section is to clarify the relationship between the
notions of formal Morita equivalence \cite{BuWa99a}, Morita
$^*$-equivalence \cite{Ara99a} and strong Morita equivalence 
\cite{Rief82,RaeWil98}. We start by briefly recalling some of the
definitions.

Let $\mathcal A$ be a $C^*$-algebra. A \emph{right pre-Hilbert 
$\mathcal A$-module}
is a right $\mathcal A$-module $\mathfrak X$, 
equipped with a positive definite $\mathcal A$-valued inner product
$\SPA{\cdot,\cdot}$
We call $\mathfrak{X}_{\mathcal A}$ a 
\emph{right Hilbert $\mathcal A$-module}
if $\mathfrak X$ is in addition complete under the norm 
$\|x\|_{\mathcal A} = \|\SPA{x,x}^{1/2} \|$. A 
\emph{left Hilbert $\mathcal A$-module} is defined analogously.
A Hilbert $\mathcal A$-module is called \emph{(topologically) full} if
the ideal $\SPA{\mathfrak X,\mathfrak X}$ is dense in $\mathcal A$.
If $\mathcal A$ and $\mathcal B$ are $C^*$-algebras, then an
\emph{imprimitivity bimodule} is a 
$(\mathcal{B}$-$\mathcal{A})$-bimodule $\BXA$
which is a full right $\mathcal A$-Hilbert module and a full left
$\mathcal{B}$-Hilbert module so that for all $x,y,z \in \BXA$ and 
$A \in \mathcal A$, $B \in \mathcal B$, we have $\SPA{B\cdot x,y}=
\SPA{x, B^*\cdot y}$, $\BSP{x\cdot A,y} = \BSP{x, y \cdot A^*}$ and
$x\cdot \SPA{y,z} = \BSP{x,y}\cdot z$.
We finally say that two $C^*$-algebras $\mathcal A$ and $\mathcal B$ are
\emph{strongly Morita equivalent} if there is a 
$(\mathcal{B}$-$\mathcal{A})$-imprimitivity bimodule. An equivalent
definition of strong Morita equivalence can be given as follows:
$\mathcal A$ and $\mathcal B$ are strongly Morita equivalent if there
exists a full right $\mathcal A$-Hilbert module 
$\mathfrak{X}_{\mathcal A}$ so that 
$\mathcal{B} \cong \mathcal{K}(\mathfrak{X}_{\mathcal A})$, where
$\mathcal{K}(\mathfrak{X}_{\mathcal A})$ is the completion of the
$^*$-algebra $\mathfrak{F}(\mathfrak{X}_{\mathcal A})$ of
``finite-rank operators'' on $\mathfrak{X}_{\mathcal A}$, that is, the
linear span of operators of the form $\Theta_{x,y}$, 
where $\Theta_{x,y} \cdot z=x \cdot \SPA{y,z}$,
for $x,y,z \in \mathfrak{X}_{\mathcal A}$.

To state Ara's result on the algebraic characterization of strong
Morita 
equivalence, we will need the following definitions. 
Recall that a ring $R$  is called \emph{non-degenerate} if 
$r \cdot R=0$ or 
$R\cdot r = 0$ implies that $r=0$ for all $r \in R$, and
\emph{idempotent} if elements of the form $r_1 r_2$ span $R$.
We also consider the category $\mathsf{Mod}$-$R$ of right 
$R$-modules $\mathcal E$ satisfying $\mathcal{E}\cdot R =
\mathcal{E}$ (in our terminology, this means that the action of $R$ on
$\mathcal E$ is strongly non-degenerate) and $x\cdot R=0$ implies
$x=0$ for all $x\in \mathcal E$, i.e. the $R$-action on $\mathcal E$
is non-degenerate. Let now $R$ and $S$ be two 
non-degenerate and idempotent rings with involution.
Then $R$ and  $S$ are called \emph{Morita $^*$-equivalent} if there
exists a right $R$-module $\mathcal E \in \mathsf{Mod}$-$R$ equipped
with a full $R$-valued inner product
$\SP{\cdot,\cdot}$ so that $\SP{x,\mathcal{E}}=0$ implies $x=0$ for 
$x,y \in \mathcal E$ and $  S \cong \mathfrak{F}(\mathcal{E}_{  R})$.
Here $\mathfrak{F}(\mathcal{E}_{  R})$ again denotes the $^*$-algebra of
``finite rank'' operators on $\mathcal{E}_{  R}$ and is defined as before
for $C^*$-algebras. We remark that if $R$ is a $K$-algebra, for
some fixed unital and 
commutative ring $K$, and $\mathcal{E}$ has a compatible $K$-module
structure, then it follows from the non-degeneracy of $R$ that 
$\SP{\cdot,\cdot}$ is automatically $K$-linear in the second
argument.

With this notion of Morita $^*$-equivalence, one can give a purely
algebraic characterization of strong Morita equivalence of
$C^*$-algebras. Recall that if $\mathcal A$ is a $C^*$-algebra, then
its \emph{Pedersen ideal} $\mathcal{P}_{\mathcal A}$  is defined as
the minimal dense ideal of $\mathcal{A}$, see
\cite[Sect.~5.6]{Peder79}. One has then the following result, due to
Ara (\cite[Thm.~2.4]{Ara99b}): 
\emph{Two $C^*$-algebras $\mathcal A$ and $\mathcal B$ are strongly
  Morita equivalent if and only if $\mathcal{P}_\mathcal{A}$ and
  $\mathcal{P}_\mathcal{B}$ are Morita $^*$-equivalent}. 
We will now discuss how this result is related to the 
notion of formal Morita equivalence. We start with the following
observation. 
\begin{lemma}\label{FormalImpliesStarLem}
Let $\mathcal A$ and $\mathcal B$ be non-degenerate, idempotent
$^*$-algebras over $\ring C$. Then if $\mathcal A$ and $\mathcal B$ are 
formally Morita equivalent, they are also Morita $^*$-equivalent.
\end{lemma}
\begin{proof}
Note that it suffices to show that there exists an equivalence
bimodule $\BXA$ so that $\SPA{x, \mathfrak{X}}=0$ implies $x=0$, for 
$x\in \mathfrak{X}$, 
$(\mathfrak{X}_{\mathcal{A}}, \SPA{\cdot,\cdot}) 
\in \mathsf{Mod}$-$\mathcal A$ 
(i.e., the $\mathcal A$-action on
$\mathfrak{X}$ is strongly non-degenerate and 
$x \cdot \mathcal{A}=0$ implies 
$x = 0$ for $x \in \mathfrak{X}$) and 
$\mathcal B \cong \mathfrak{F}(\mathfrak{X}_{\mathcal A})$.

Let $\BXA$ be an arbitrary equivalence bimodule.
We start by noticing that we can choose it
so that the $\mathcal B$ and $\mathcal A$ actions are strongly
non-degenerate. Indeed, it follows from the compatibility of 
$\SPA{\cdot,\cdot}$ and $\BSP{\cdot,\cdot}$ and their fullness
properties that 
$\mathfrak{X}\cdot \mathcal{A} = \mathcal{B}\cdot \mathfrak{X}$. If we
denote $\hat{\mathfrak X}=
\mathfrak{X}\cdot \mathcal{A} = \mathcal{B} \cdot \mathfrak{X}$, then
$\hat{\mathfrak X}$ has a natural $(\mathcal B$-$\mathcal A)$-bimodule
structure and also compatible $\mathcal A$ and $\mathcal B$ valued
inner products, defined by simply restricting $\BSP{\cdot,\cdot}$ and 
$\SPA{\cdot,\cdot}$ to $\hat{\mathfrak X}$. Note that given $A \in
\mathcal A$, by idempotency we can write $A = \sum_i A^1_iA^2_iA^3_i$,
for $A^j_i \in \mathcal A$, and by fullness of $\SPA{\cdot,\cdot}$, we
can write $A^2_i=\sum_j \SPA{x^i_j,y^i_j}$. So  
$A = \sum_i A^1_i\sum_j \SPA{x^i_j,y^i_j} A^3_i = \sum_{i,j}
\SPA{x^i_j  {A^1_i}^*, y^i_j  A^3_i}$, showing that 
$\SPA{\cdot,\cdot}$ restricted to $\hat{\mathfrak X}$ is still full, 
and the same holds for $\SPB{\cdot,\cdot}$. Therefore 
$\hat{\mathfrak X}$ is a $(\mathcal B$-$\mathcal A)$-equivalence
bimodule with strongly non-degenerate $\mathcal A$ and $\mathcal B$
actions.

We also observe that if $\BXA$ is an (arbitrary) equivalence bimodule,
then $\mathcal B$ and $\mathcal A$ act on $\mathfrak{X}$ faithfully
(i.e., $B \cdot x =0$ for all $x \in \mathfrak{X}$ implies $B=0$, for
all $B \in \mathcal B$, and the same for $\mathcal A$). To see that,
suppose $B \cdot x =0$ for all $x$. Then given 
an arbitrary $B' \in \mathcal B$, it follows that we can write 
$B'=\sum_i\BSP{x_i,y_i}$ and hence $B\cdot B' = B \sum_i\BSP{x_i,y_i}
= \sum_i\BSP{B x_i,y_i} =0$. Therefore $B=0$ by non-degeneracy of
$\mathcal B$ and the same argument holds for $\mathcal A$. 
It then follows (see the proof of \cite[Prop.~5.16]{BuWa99a}) that
$\{x \in \BXA | \SPA{x,y}=0, \, \forall y \in \BXA \} =
\{x \in \BXA | \BSP{x,y}=0, \, \forall y \in \BXA \} = N$ and
$\mathfrak{X}/N$ is still a $(\mathcal B$-$\mathcal A)$-equivalence
bimodule, with the property that $\SPA{x, \mathfrak{X}}=0$ implies
that $x=0$, the same holding for $\SPA{\cdot,\cdot}$. Moreover, if the
actions of $\mathcal A$ and $\mathcal B$ were strongly non-degenerate
on $\BXA$, this still holds for the quotient $\mathfrak{X}/N$. So, we
can always find an equivalence bimodule $\BXA$ with strongly
non-degenerate $\mathcal A$ and $\mathcal B$ actions and inner
products satisfying $\SPA{x, \mathfrak{X}}=0 \implies x=0$ and 
$\BSP{x, \mathfrak{X}}=0 \implies x=0$. We observe that, in this case 
(see proof of \cite[Prop.~6.1]{BuWa99a}), $\mathcal B \cong 
\mathfrak{F}(\mathfrak{X}_{\mathcal A})$.

Given $\BXA$ an equivalence bimodule with these properties, we note
that if $x \cdot \mathcal{A}=0$, then in 
particular $x\cdot \SPA{y,z} = 0$ for all  
$y,z \in \mathfrak{X}$. But then $\BSP{x,y}\cdot z =0$ for all $z$,
and hence $\BSP{x,y}=0$. But since $y$ is also arbitrary, it follows
that $x=0$. Therefore, 
$(\mathfrak{X}_{\mathcal{A}},\SPA{\cdot,\cdot}) \in 
\mathsf{Mod}$-$\mathcal A$ and clearly establishes a Morita
$^*$-equivalence between $\mathcal A$ and 
$\mathfrak{F}(\mathfrak{X}_{\mathcal A}) \cong \mathcal B$.
\end{proof}

We observe that $^*$-algebras with approximate identities are
non-degenerate and idempotent. We will now show the converse of the
previous lemma in the case of $C^*$-algebras and Pedersen ideals
(recalling that both types of algebras are also non-degenerate and
idempotent). We start with some general lemmas. 
\begin{lemma} \label{PositivityLem}
Let $\mathcal A$ be a $C^*$-algebra and $\mathcal{P}_{\mathcal A}$ its
Pedersen ideal. Then $A \in \mathcal{P}_{\mathcal A}$ is positive in
$\mathcal A$ if and only if $A \in {\mathcal{P}_{\mathcal A}}^+$
(positivity in the purely algebraic sense of $^*$-algebras over 
$\ring C$). 
\end{lemma}
\begin{proof}
Let $A \in \mathcal{P}_{\mathcal A}$ be positive in $\mathcal A$ and let
$\omega: \mathcal{P}_{\mathcal A} \longrightarrow \mathbb{C}$ be an
arbitrary positive linear functional (in the algebraic sense, that is,
$\omega (A^*A)\geq 0$ for all $A \in \mathcal{P}_{\mathcal A}$). We must
show that $\omega (A) \geq 0$. Note that by the general properties of
Pedersen ideals (see \cite[Thm.~5.6.2]{Peder79}), if 
$A \in \mathcal{P}_{\mathcal A}$, then 
$C^*(A) \subseteq \mathcal{P}_{\mathcal A}$, where $C^*(A)$ denotes
the $C^*$-algebra generated by $A$. Also observe that if 
$A \in \mathcal{A}^+$, then $A \in C^*(A)^+$, and hence $A=A_1^*A_1$,
for some $A_1 \in C^*(A) \subseteq \mathcal{P}_{\mathcal A}$.
Therefore, $\omega(A)=\omega(A_1^*A_1) \geq 0$.
For the converse, note that if 
$\omega$ is a positive linear functional in $\mathcal A$, then 
$\omega|_{\mathcal{P}_{\mathcal A}}$ is positive in
$\mathcal{P}_{\mathcal A}$ and hence if 
$A \in \mathcal{P}_{\mathcal A}^+$, we immediately have 
$\omega(A) \geq 0$. Hence $A \in \mathcal A^+$.
\end{proof}

We now observe a general property of $^*$-representations of Pedersen
ideals on pre-Hilbert spaces which should be well-known.
\begin{lemma}\label{ExtensionLem}
Let $\mathcal A$ be a $C^*$-algebra,
$\mathcal{P}_{\mathcal A}$ its Pedersen ideal and suppose
$\pi:\mathcal{P}_{\mathcal A} \longrightarrow \Bounded(\mathfrak{H})$
is a $^*$-representation on a complex pre-Hilbert space
$\mathfrak{H}$.
Then 
$\pi$  extends to a $^*$-representation 
$\pi^\cl: \mathcal{P}_{\mathcal A} \longrightarrow \Bounded(\mathfrak H^\cl)$,
where $\mathfrak H^\cl$ is the closure of $\mathfrak H$. Moreover,
$\|\pi^\cl(A)\| \leq \|A\|$ for all $A \in \mathcal{P}_{\mathcal A}$ and hence
$\pi^\cl$ also extends to a $^*$-representation of $\mathcal A$ on 
$\mathfrak H^\cl$.
\end{lemma} 
\begin{remark}
As pointed out to us by P. Ara, the above lemma follows from a more general 
fact. If $\mathcal A$ is any complex $^*$-algebra with positive definite
involution (\cite[pp.~338]{Hand81}), let $\mathcal{A}_b$
be the $^*$-subalgebra of bounded elements of $\mathcal A$ 
(\cite[pp.~338]{Hand81}). Then any $^*$-representation of $\mathcal{A}_b$
on a complex pre-Hilbert space $\mathfrak H$ naturally extends to a 
$^*$-representation of $\overline{\mathcal{A}_b}$ (closure of $\mathcal{A}_b$
with respect to its natural seminorm, as defined in \cite[pp.~342]{Hand81})
on $\mathfrak H^\cl$.
\end{remark}

We can now discuss general positivity properties of modules over
Pedersen ideals.
\begin{proposition} \label{PropertyPProp}
Let $\mathcal A$ be a $C^*$-algebra, $\mathcal{P}_{\mathcal A}$ its Pedersen
ideal and $\mathcal B$ an arbitrary $^*$-algebra over $\mathbb{C}$.
Suppose $_\mathcal{B}\mathfrak{X}_{\mathcal{P}_{\mathcal A}}$ is a bimodule,
equipped with a positive semi-definite $\mathcal{P}_{\mathcal A}$-valued
inner product $\SP{\cdot, \cdot}_{\mathcal{P}_{\mathcal A}}$
and compatible $\mathcal B$-action. Then $\mathfrak{X}$ satisfies
property \axiom{P}.
\end{proposition}
\begin{proof}
We must show that given a complex pre-Hilbert space $\mathfrak H$
and a $^*$-representation $\pi: \mathcal{P}_{\mathcal A}
\longrightarrow \Bounded(\mathfrak H)$,
the formula $\SPKT{x_1\otimes \phi_1,x_2\otimes \phi_2}
= \SPH{\phi_1,\pi(\SP{x_1, x_2}_{\mathcal{P}_{\mathcal A}}) \phi_2}$,
for $x_1,x_2 \in \mathfrak{X}$ and $\phi_1,\phi_2 \in \mathfrak H$ defines
a positive semi-definite inner product on $\tildeK = 
\mathfrak{X} \otimes_{\mathcal{P}_{\mathcal A}} \mathfrak{H}$. 
Note that since $\pi$ extends to a $^*$-representation $\pi^\cl:
\mathcal{P}_{\mathcal A}\longrightarrow \Bounded(\mathfrak H^\cl)$ on the
Hilbert space completion of $\mathfrak H$, there is a natural isometric
embedding $\iota: \mathfrak{X}\otimes_{\mathcal{P}_{\mathcal A}}\mathfrak H
\longrightarrow \mathfrak{X}\otimes_{\mathcal{P}_{\mathcal A}} 
\mathfrak H^\cl$, where the $\mathcal{P}_{\mathcal A}$-balanced tensor products
are defined using $\pi$ and $\pi^\cl$, respectively. So for
$\phi_1,\ldots,\phi_n \in \mathfrak H \subseteq \mathfrak H^\cl$ and
$x_1,\ldots,x_n \in \mathfrak{X}$, we have 
\begin{equation}\label{IsometricEq}
    \SPKT{\sum_i x_i\otimes\phi_i,\sum_i x_i\otimes\phi_i} 
    = 
    \sum_{i,j} \SPH{\phi_i, \pi(\SPA{x_i, x_j})\phi_j}
    = 
    \sum_{i,j} 
    \SP{\phi_i, \pi^\cl(\SPA{x_i, x_j})\phi_j}_{\mathfrak H^\cl}. 
\end{equation}
Therefore, to show the positivity of $\SPKT{\cdot,\cdot}$, it suffices
to show that the right hand side of (\ref{IsometricEq}) is
non-negative. This can  be done by the usual procedure (see
\cite[Ch.~IV, Sect.~2.2]{Land98}), 
using the fact that any $^*$-representation of
$\mathcal{P}_{\mathcal A}$ on a Hilbert space can be decomposed into
the sum of a trivial representation and a non-degenerate one, 
and the non-degenerate part itself can be decomposed into the sum 
(in the topological sense) of cyclic (also in the topological sense)
$^*$-representations (see \cite[Ch.~I, Prop.~1.5.2]{Land98}). 
\end{proof}

Let $\mathcal A$ and $\mathcal B$ be $C^*$-algebras,
let $\mathcal{P}_{\mathcal A}$ and $\mathcal{P}_{\mathcal B}$ be the 
corresponding Pedersen ideals and suppose there exists 
$\mathfrak{X} \in \mathsf{Mod}$-$\mathcal{P}_{\mathcal A}$, 
equipped with a pairing $\SP{\cdot,\cdot}$
so that $(\mathfrak{X}, \SP{\cdot,\cdot})$ establishes a Morita
$^*$-equivalence between $\mathcal{P}_{\mathcal A}$ and 
$\mathfrak{F}(\mathfrak{X}_{\mathcal{P}_{\mathcal A}}) \cong 
\mathcal{P}_{\mathcal B}$. In \cite[Thm.~2.4]{Ara99b}, it is shown
that, in this case, one automatically has 
$\SP{x,x} \in \mathcal{A}^+$  and $\Theta_{x,x} \in \mathcal{B}^+$
(recall that for $C^*$-algebras, the usual notion of positivity
coincides with  positivity as defined for arbitrary $^*$-algebras 
over $\ring C$). Then, by Lemma \ref{PositivityLem}, it follows that
$\SP{x,x} \in \mathcal{P}_{\mathcal A}^+$ and $\Theta_{x,x}
\in \mathcal{P}_{\mathcal B}^+$. It is then easy to check (see
e.g. \cite[Prop.~6.2]{BuWa99a}) that $\mathfrak{X}$ is a 
$(\mathcal{P}_{\mathcal B}$-$\mathcal{P}_{\mathcal A})$-bimodule
satisfying all the properties of an equivalence bimodule, except
possibly for \axiom{P} and \axiom{Q}. But these properties follow
directly from Proposition \ref{PropertyPProp}. Hence we have
\begin{corollary}
If $\mathcal{P}_{\mathcal A}$ and $\mathcal{P}_{\mathcal B}$ are
Morita $^*$-equivalent, then they are automatically formally Morita
equivalent. 
\end{corollary}
It then follows from Lemma \ref{FormalImpliesStarLem},
Proposition \ref{PropertyPProp} and \cite[Thm.~2.4]{Ara99b} that
\begin{theorem}
\label{TheNiceTheoremAboutCStar}
Two $C^*$-algebras $\mathcal A$ and $\mathcal B$ are strongly Morita
equivalent if and only if their corresponding Pedersen ideals,
$\mathcal{P}_{\mathcal A}$ and $\mathcal{P}_{\mathcal B}$, are
formally Morita equivalent. In particular, if $\mathcal A$ and
$\mathcal B$ are unital, then they are strongly Morita equivalent if
and only if they are formally Morita equivalent.
\end{theorem}

It follows from \cite[Thm.~5.10]{BuWa99a} that if
two $C^*$-algebras $\mathcal A$ and $\mathcal B$
are strongly Morita equivalent, then 
$^*$-$\mathsf{Rep}(\mathcal{P}_{\mathcal A})$
and $^*$-$\mathsf{Rep}(\mathcal{P}_{\mathcal B})$ are equivalent
categories. We recall from \cite{Rief74b} 
that strong Morita equivalent
$C^*$-algebras have equivalent categories of Hermitian modules. We conclude
this section discussing how these two results are related.

For a $C^*$-algebra $\mathcal A$, let $\mathsf{hermod}(\mathcal A)$ 
be the category of $^*$-representations of $\mathcal A$ on Hilbert 
spaces, with isometric intertwiners as morphisms. The category of
\emph{non-degenerate} $^*$-representations is denoted by 
$\mathsf{Hermod}(\mathcal A)$ and objects in this category are called
\emph{Hermitian modules} \cite{Rief74}. In order to understand how these
categories are related to the categories of purely algebraic representations
$^*$-$\mathsf{rep}(\mathcal{P}_{\mathcal A})$ and 
$^*$-$\mathsf{Rep}(\mathcal{P}_{\mathcal A})$, we introduce two other
categories. By $\srepDSA$ we denote the
category of $^*$-representations of $\mathcal A$ on Hilbert spaces with a 
dense $\mathcal{P}_{\mathcal A}$-invariant subspace, while for $\sRepDSA$ 
we require in addition that the 
induced $^*$-representation of $\mathcal{P}_{\mathcal A}$
restricted to the dense subspace be strongly non-degenerate
(in the algebraic sense). For the morphisms we take isometric intertwiners
preserving the dense invariant subspaces. Note that objects in $\sRepDSA$ are
in particular Hermitian modules over $\mathcal A$ 
(with a distinguished dense subspace).
 One can check that the process of
extending the $^*$-representations, discussed in Lemma
\ref{ExtensionLem}, gives rise to a functor from 
$^*$-$\mathsf{Rep}(\mathcal{P}_{\mathcal A})$ to $\sRepDSA$ 
($^*$-$\mathsf{rep}(\mathcal{P}_{\mathcal A})$ to $\srepDSA$) for which the 
natural restriction functor is an inverse. It then follows that 
$^*$-$\mathsf{Rep}(\mathcal{P}_{\mathcal A})$ and $\sRepDSA$ 
($^*$-$\mathsf{rep}(\mathcal{P}_{\mathcal A})$ and $\srepDSA$)
are equivalent categories. We point out that $\sRepDSA$ and 
$\mathsf{Hermod}(\mathcal A)$ ($\srepDSA$ and 
$\mathsf{hermod}(\mathcal A)$) are not equivalent categories in general.
We finally remark that the categories $\srepDSA$ and $\sRepDSA$ are interesting
in their own right since,
in many applications,  $^*$-representations of $\mathcal A$ come
automatically with a dense invariant subspace, as in GNS representations. 
As another example, certain unbounded operators associated to 
$^*$-representations, such as generators of symmetries, distinguish  dense
subspaces by their domains of definition.

%
%

\section{Closed $^*$-Ideals and the Minimal Ideal}
\label{MinimalSec}

After having clarified the connection between formal Morita
equivalence and $C^*$-algebras, we consider in this section general
$^*$-algebras over $\ring C$ and start the discussion about (closed)
$^*$-ideals.

If $\mathcal I \subseteq \mathcal A$ is a $^*$-ideal, i.e. an ideal
closed under the $^*$-involution, then $\mathcal A \big/ \mathcal I$
is again a $^*$-algebra. If $\mathcal I$ is a left (or right) ideal
which is closed under the $^*$-involution then it is automatically a
two-sided ideal and hence a $^*$-ideal.
\begin{definition}
A $^*$-ideal $\mathcal I$ of $\mathcal A$ is called closed if it is
the kernel of a $^*$-representation.
\end{definition}
Clearly, the kernel $\ker\pi$ of a $^*$-representation $\pi$ is a
$^*$-ideal, and unitary equivalent $^*$-re\-pre\-sen\-ta\-tions have
the same kernel. But $\ker\pi$ does not determine this equivalence
class since e.g. $\ker\pi = \ker(\pi\oplus\pi)$. 
We observe that if $\mathcal A$ is a $C^*$-algebra, then a $^*$-ideal
$\mathcal I$ in $\mathcal A$ is closed in this algebraic sense if and
only if it is norm closed. Indeed, recall that a $^*$-ideal in a
$C^*$-algebra is norm closed if and only if it is the kernel of some
$^*$-representation of $\mathcal A$ on a Hilbert space (\cite[Ch.~I,
Thm.~1.3.10]{Land98}). But since any representation of a $C^*$-algebra
on a pre-Hilbert space extends to a representation on its completion
(see Lemma \ref{ExtensionLem}),
the conclusion follows. 
\begin{lemma}
\label{IntersectClosedIdealLem}
Let $\mathcal A$ be a $^*$-algebra over $\ring C$ and $\mathcal I$ a
$^*$-ideal. Then $\mathcal A \big/ \mathcal I$ has a faithful
$^*$-representation if and only if $\mathcal I$ is closed. If
$\{\mathcal I_\alpha\}_{\alpha \in \Lambda}$ are closed $^*$-ideals
then $\mathcal I = \bigcap_{\alpha \in \Lambda} \mathcal I_\alpha$ is
again a closed $^*$-ideal.
\end{lemma}
\begin{proof}
For the first part use the projection 
$\rho: \mathcal A \to \mathcal A \big/ \mathcal I$ to pull-back
$^*$-representations. For the second observe that if 
$\mathcal I_\alpha = \ker \pi_\alpha$ then 
$\mathcal I = \ker (\bigoplus_{\alpha \in \Lambda} \pi_\alpha)$.
\end{proof}

Since formal Morita equivalence deals with strongly non-degenerate
$^*$-representations, it would be nice if these particular
representations were sufficient to define all closed $^*$-ideals. If
$\mathcal A$ is idempotent, then this is indeed the case.  
\begin{lemma}
\label{SNDLem}
Let $\mathcal A$ be an idempotent $^*$-algebra over $\ring C$ and
$\pi: \mathcal A \to \Bounded (\mathfrak H)$ a
$^*$-representation. Let 
\begin{equation}
\label{SNDHilbertSpaceDef}
    \mathfrak H_{\snd} := 
    \big\{\phi \in \mathfrak H \; \big| \;
    \exists \phi_i \in \mathfrak H, A_i \in \mathcal A, 
    i = 1, \ldots, n:\;
    \phi = \mbox{$\sum_i$} \pi(A_i) \phi_i \big\}.
\end{equation}
Then $\mathfrak H_{\snd}$ is $\pi$-invariant and the restriction
$\pi_{\snd} := \pi|_{\mathfrak H_{\snd}}$ is strongly non-degenerate
with $\ker \pi_{\snd} = \ker\pi$. Moreover, every GNS representation
is strongly non-degenerate.
\end{lemma}
\begin{proof}
The invariance of $\mathfrak H_{\snd}$ is obvious and the fact that
$\pi_{\snd}$ is strongly non-degenerate follows from the idempotency
of $\mathcal A$. Finally let $\pi_{\snd} (A) = 0$. Then in particular
$\pi(A)\pi(A^*)\phi = 0$ for all $\phi \in \mathfrak H$ whence 
$\pi(A) = 0$ follows and thus $\ker\pi_{\snd} \subseteq \ker\pi$. The
other inclusion is trivial. The strong non-degeneracy of a GNS
representation follows directly from the definition.
\end{proof}

Using the closed $^*$-ideals, we define the following closure
operation. If $\mathcal J \subseteq \mathcal A$ is an arbitrary subset
then $\mathcal J^\cl$ is the smallest closed $^*$-ideal containing
$\mathcal J$, i.e.
\begin{equation}
\label{ClosureDef}
    \mathcal J^\cl := 
    \bigcap_{\mathcal J \subseteq \mathcal I, \:  
             \mathcal I \text{ closed $^*$-ideal}} \mathcal I,
\end{equation}
where $\mathcal I$ runs over the closed $^*$-ideals of 
$\mathcal A$. Clearly $\mathcal J \subseteq \mathcal J^\cl$ and
$(\mathcal J^\cl)^\cl = \mathcal J^\cl$ for any subset of 
$\mathcal A$. This implies that for
$\mathcal J \subseteq \mathcal I$ we have
$\mathcal J^\cl \subseteq \mathcal I^\cl$ and 
$(\bigcup_{\alpha \in \Lambda} \mathcal I_\alpha)^\cl 
= (\bigcup_{\alpha \in \Lambda} \mathcal I_\alpha^\cl)^\cl$
for an arbitrary index set $\Lambda$ and 
$\mathcal I_\alpha \subseteq \mathcal A$. 
Finally note that 
$\bigcap_{\alpha \in \Lambda} \mathcal I_\alpha^\cl 
= \left(\bigcap_{\alpha \in \Lambda} 
\mathcal I_\alpha^\cl\right)^\cl$.
We can define the following operations. For
two subsets $\mathcal I, \mathcal J \subseteq \mathcal A$, we set
\begin{equation}
\label{LatticeDef}
    \mathcal I \vee \mathcal J := (\mathcal I \cup \mathcal J)^\cl
    \quad
    \text{ and }
    \quad
    \mathcal I \wedge \mathcal J := \mathcal I^\cl \cap \mathcal J^\cl.
\end{equation}
One can check that these two operations define a lattice structure on
the set of closed $^*$-ideals of $\mathcal A$. In this case,
the operation `$\le$', which is defined by $\mathcal I \le \mathcal J$ if
$\mathcal I \wedge \mathcal J = \mathcal I$, 
coincides with the set-theoretic `$\subseteq$'.
We then have a lattice structure on the set
of all closed $^*$-ideals in $\mathcal A$, that we denote by
$\mathfrak{L}(\mathcal A)$. Note that all $^*$-ideals,
closed or not, also have a lattice structure, with the `closure map'
defined by taking the smallest $^*$-ideal containing the corresponding
subset. However, we will be mainly interested in the closed
$^*$-ideals.

For a $C^*$-algebra $\mathcal A$, $\mathfrak{L}(\mathcal A)$ coincides
with its usual lattice of norm closed $^*$-ideals. If
$\mathcal{P}_{\mathcal A}$ is its Pedersen ideal, note that we have a
natural map  
$F:\mathfrak{L}(\mathcal{P}_{\mathcal A}) 
\to \mathfrak{L}(\mathcal A)$ 
defined as follows. Let 
$\mathcal I \in \mathfrak{L}(\mathcal{P}_{\mathcal A})$. Then
$\mathcal I=\ker \pi$, for some $^*$-representation of
$\mathcal{P}_{\mathcal A}$ on a pre-Hilbert space $\mathfrak H$. As we
observed in Lemma \ref{ExtensionLem}, $\pi$ extends to a
$^*$-representation $\pi^\cl$ of $\mathcal A$ on  
$\mathfrak H^\cl$. We then define $F$ as the map 
$\mathcal I=\ker \pi \mapsto F(\mathcal I)= \ker \pi^\cl$.
\begin{proposition}\label{SameLatticeProp}
$F:\mathfrak{L}(\mathcal{P}_{\mathcal A}) \longrightarrow
\mathfrak{L}(\mathcal A)$ is a lattice isomorphism. 
\end{proposition}
\begin{proof}
Let us first observe that $F$ is a well-defined map. 
Note that if $\pi$ is a $^*$-representation of 
$\mathcal{P}_{\mathcal A}$, then 
$\ker\pi$ is always norm dense in $\ker \pi^\cl$. To see that, let 
$\{e_\lambda\}_{\lambda \in \Lambda}$ be an approximate identity in 
$\mathcal A$, with $e_\lambda \in \mathcal{P}_{\mathcal A}$. 
This is possible since $\mathcal{P}_{\mathcal A}$ is dense in 
$\mathcal A$ (\cite[Thm.~I.7.2]{Dixm77}). Then given 
$A \in \ker \pi^\cl$, we have 
$e_\lambda A \in \ker\pi^\cl \bigcap \mathcal{P}_{\mathcal A} = \ker\pi$ 
and $e_\lambda A \to A$. Now let $\mathcal I \in 
\mathfrak{L}(\mathcal{P}_{\mathcal A})$ and suppose
$\pi$ and $\rho$ are two different $^*$-representations of 
$\mathcal{P}_{\mathcal A}$ so that $\mathcal I = \ker \pi = \ker \rho$.
We must show that $\ker \pi^\cl=\ker \rho^\cl$. 
Note that if $A \in \ker \pi^\cl$, then we can find $A_\lambda \in \ker \pi
=\ker \rho$ with $A_\lambda \to A$. So $\rho(A_\lambda)=0$, for all
$\lambda$, and hence $\rho(A)=0$. This shows that $\ker \pi^\cl
\subseteq \ker \rho^\cl$ and by changing the roles of $\pi$ and $\rho$
we conclude that $\ker \pi^\cl = \ker \rho^\cl$. This proves
well-definedness. It is not hard to check that $F$ is a bijection
which preserves the order structure of the lattices. This guarantees
that $F$ is a lattice isomorphism. 
\end{proof}

Recall that for a $C^*$-algebra, any (topologically) non-degenerate 
$^*$-representation can be decomposed into a direct sum (also in the
topological sense) of GNS representations. Thus for $C^*$-algebras,
any closed $^*$-ideal can be obtained as the kernel of a (topological)
sum of GNS representations. 
We observe that for a general $^*$-algebra over $\ring C$, it is not
true that an arbitrary (strongly non-degenerate) $^*$-representation
decomposes into an algebraic sum of GNS representations. In fact, even
if $\ring C$ is an algebraically closed field (which is complete in
its canonical topology) and we allow topological sums, this decomposition
does not hold in general. This follows since otherwise every Hilbert 
space over $\ring C$ would have a Hilbert basis (this notion makes
sense in this situation), which is \emph{not} the case (see
e.g.~\cite[Thm.~7]{BW98a}). Nevertheless, it is still true that
arbitrary closed $^*$-ideals on a $^*$-algebra arise as kernel of sums
of GNS representations, as we now observe.

Let $\omega: \mathcal A \to \ring C$ be a positive linear
functional. For $B \in \mathcal A$, we consider the positive linear
functional 
$\omega_B (A) := \omega(B^*AB) = \SP{\psi_B, \pi_\omega(A)\psi_B}$.
Then $A \in \ker\pi_\omega$ if and only if $\omega_B(A) = 0$ for all
$B \in \mathcal A$.  Using a
polarization argument and the fact that 
$\Bounded(\mathfrak H_\omega)$ is torsion-free (\cite[Prop.~2.8]{BuWa99a}),
we see that this is equivalent to 
$\omega_B(A^*A) = 0$ for all $B \in \mathcal A$. Hence we have 
\begin{equation}
\label{KerGNSRep}
    \ker \pi_\omega 
    = \bigcap_{B \in \mathcal A} \mathcal J_{\omega_B}
    = \bigcap_{B \in \mathcal A} \ker\omega_B.
\end{equation}
Note that neither $\mathcal J_{\omega_B}$ nor $\ker\omega_B$ is a
$^*$-ideal in general but only the above intersection. Now consider an
arbitrary $^*$-representation $\pi$ on $\mathfrak H$. For every
$\varphi \in \mathfrak H$, one defines the expectation value
$\omega_\varphi(A) := \SP{\varphi, \pi(A)\varphi}$, which is clearly  a
positive linear functional of $\mathcal A$. Denoting the Gel'fand
ideal by $\mathcal J_\varphi$ and the corresponding GNS representation
by $\pi_\varphi$, it follows that $A \in \ker\pi_\varphi$ if and only
if $\pi(A)\pi(B)\varphi = 0$ for all $B \in \mathcal A$. But
this implies that $A \in \mathcal J_\varphi$ by taking $B = A^*$. Moreover,
$A \in \mathcal J_\varphi$ if and only if $\pi(A)\varphi = 0$. So
$A \in \ker\omega_\varphi$. Thus we have the inclusions
\begin{equation}
\label{ThreeSpecialInclusions}
    \ker\pi_\varphi 
    \subseteq \mathcal J_\varphi 
    \subseteq \ker\omega_\varphi.
\end{equation}
\begin{lemma}
\label{KernelGNSKernelLem}
Let $\mathcal A$ be a $^*$-algebra over $\ring C$ and 
$\pi: \mathcal A \to \Bounded(\mathfrak H)$ a
$^*$-representation. Consider 
$\pi_M := \bigoplus_{\varphi \in \mathfrak H} \pi_\varphi$.
Then one has 
\begin{equation}
\label{KerPiKerPiMonster}
    \ker\pi 
    = \ker\pi_M 
    = \bigcap_{\varphi \in \mathfrak H} \ker\pi_\varphi
    = \bigcap_{\varphi \in \mathfrak H} \mathcal J_\varphi
    = \bigcap_{\varphi \in \mathfrak H} \ker \omega_\varphi.
\end{equation}
\end{lemma}
\begin{proof}
It remains to show 
$\bigcap_{\varphi \in \mathfrak H} \ker \omega_\varphi = \ker\pi$.
But this follows by polarization and the fact that 
$\Bounded(\mathfrak H)$ is torsion-free.
\end{proof}
\begin{lemma}
Let $\mathcal A$ be a $^*$-algebra over $\ring C$ with an approximate
identity. Then for every positive linear functional $\omega$ one has
\begin{equation}
\label{ThreeInclusion}
    \ker\pi_\omega 
    \subseteq \mathcal J_\omega
    \subseteq \ker \omega.
\end{equation}
\end{lemma}
\begin{proof}
This generalizes (\ref{ThreeSpecialInclusions}). Let 
$A \in \mathcal A$ and $E_\alpha = E^*_\alpha \in \mathcal A$ such
that $E_\alpha A = A = A E_\alpha$. Then $\pi_\omega(A) = 0$ implies
$0 = \pi_\omega (A) \psi_{E_\alpha} = \psi_A$, whence 
$A \in \mathcal J_\omega$. For $A \in \mathcal J_\omega$, we find by the 
Cauchy-Schwarz inequality that
$0 \le \omega(A) \cc{\omega(A)} 
\le \omega(A^*A)\omega(E^*_\alpha E_\alpha) = 0$ and hence
$A \in \ker\omega$.
\end{proof}

Let us investigate the lattice of closed $^*$-ideals more
closely. While the lattice of $^*$-ideals trivially has a minimal and
a maximal element, namely $\{0\}$ and $\mathcal A$, the closed
$^*$-ideals provide a more interesting structure. We define the
minimal closed $^*$-ideal of $\mathcal A$ by
\begin{equation}
\label{MinimalIdealDef}
    \Jmin (\mathcal A) := 
    \bigcap_{\mathcal I \text{ closed $^*$-ideal}}
    \mathcal I,
\end{equation}
which is again a closed $^*$-ideal due to
Lemma~\ref{IntersectClosedIdealLem}. Furthermore $\Jmin (\mathcal A)$
is clearly the minimal element of the lattice. It is obvious
that $\Jmin (\mathcal A) = \{ 0\}$ if and only if
$\mathcal A$ has a faithful $^*$-representation on some pre-Hilbert
space. Thus for a $C^*$-algebra $\mathcal A$ one automatically has
$\Jmin(\mathcal A) = \{0\}$. We note that the same holds for 
star-product algebras (see e.g.~\cite[Sect.~4]{Wal2000a}).

Since the kernel of a $^*$-representation can always be obtained
by using direct sums of GNS representations, it is sufficient to
consider GNS representations in order to characterize 
$\Jmin(\mathcal A)$. With (\ref{KerGNSRep}), this gives
\begin{equation}
\label{JminIntersections}
    \Jmin(\mathcal A)
    = \bigcap_{\omega} \ker\pi_\omega
    = \bigcap_{\omega} \bigcap_{B \in \mathcal A} \mathcal J_{\omega_B}
    = \bigcap_{\omega} \bigcap_{B \in \mathcal A} \ker\omega_B,
\end{equation}
where $\omega$ runs over all positive linear functionals. If 
$\mathcal A$ has in addition an approximate identity, we obtain the
following result:
\begin{theorem}
\label{JminTheo}
Let $\mathcal A$ be a $^*$-algebra over $\ring C$ with approximate
identity. Then the minimal closed $^*$-ideal is given by
\begin{equation}
\label{JminThreeIntersection}
    \Jmin(\mathcal A)
    = \bigcap_{\omega} \ker\pi_\omega
    = \bigcap_{\omega} \mathcal J_\omega
    = \bigcap_{\omega} \ker\omega,
\end{equation}
where $\omega$ runs over all positive linear functionals.
If $A^*A = 0$, or if $A$ is normal and nilpotent, or if $zA = 0$ for
$0 \ne z \in \ring C$ then $A \in \Jmin(\mathcal A)$. 
\end{theorem}
\begin{proof}
Let $A \in \bigcap_{\omega} \ker\omega$. Then $\pi(A) = 0$ for all
$^*$-representations $\pi$ since otherwise there would exist a
$\varphi \in \mathfrak H$ such that 
$\omega_\varphi(A) = \SP{\varphi, \pi(A)\varphi} \ne 0$ which
contradicts the assumption since $\omega_\varphi$ is a positive
functional. Hence $A \in \Jmin(\mathcal A)$ and by
(\ref{ThreeInclusion}) and (\ref{JminIntersections}) the equality
(\ref{JminThreeIntersection}) follows. The other statements follow
directly from \cite[Prop.~2.8]{BuWa99a}.
\end{proof}

Note  that all `non-$C^*$-algebra-like' elements are absorbed into
$\Jmin(\mathcal A)$. In particular
$\mathcal A \big/ \Jmin(\mathcal A)$ has sufficiently many positive
linear functionals and still an approximate identity. This can be seen
as an analogue of the construction of the $C^*$-enveloping
algebra of a Banach $^*$-algebra, see e.g.~\cite[Sect.~II.7]{Dixm77}.
The passage from $\mathcal A$ to $\mathcal A \big/ \Jmin(\mathcal A)$
is functorial in the category of $^*$-algebras with approximate
identities.
\begin{proposition}
\label{JminQuotientProp}
Let $\mathcal A$, $\mathcal B$ be $^*$-algebras over $\ring C$ with
approximate identities and let $\Phi: \mathcal A \to \mathcal B$ be a
$^*$-homomorphism. Then 
$\Phi(\Jmin(\mathcal A)) \subseteq \Jmin(\mathcal B)$ and thus there
exists a unique $^*$-homomorphism 
$\phi: \mathcal A \big/ \Jmin(\mathcal A) 
\to \mathcal B \big/ \Jmin(\mathcal B)$ such that 
\begin{equation}
\label{OurFirstCommutingDiagram}
    \begin{array} {ccc}
        \mathcal A 
        & \stackrel{\Phi}{\longrightarrow} 
        & \mathcal B \\
        \downarrow & & \downarrow \\
        \mathcal A \big/ \Jmin(\mathcal A) 
        & \stackrel{\phi}{\longrightarrow} 
        & \mathcal B \big/ \Jmin(\mathcal B) 
    \end{array}
\end{equation}
commutes. Thus the passage from $\mathcal A$ to 
$\mathcal A \big/ \Jmin(\mathcal A)$ is functorial.
\end{proposition}
\begin{proof}
If $\omega: \mathcal B \to \ring C$ is a positive linear functional
then $\Phi^*\omega = \omega \circ \Phi: \mathcal A \to \ring C$ is
positive, too, whence 
$\Phi(\Jmin(\mathcal A)) \subseteq \Jmin(\mathcal B)$ follows from
(\ref{JminThreeIntersection}). Then (\ref{OurFirstCommutingDiagram})
is clear.
\end{proof}
\begin{remark}
The minimal ideal as discussed here is related to the work of Handelman 
\cite{Hand81}, as P. Ara brought to our attention.
Let $\mathcal A$ be a (unital) $^*$-algebra over $\ring C=\ring R(\im)$ with
positive definite involution (\cite[pp.~338]{Hand81}) and suppose
$\mathbb{Q}\subseteq \ring R \subseteq \mathbb R$ (so that $\ring R$
is archimedean-ordered). 
Let $J^*(\mathcal A)$ be the set of
bounded elements of $\mathcal A$ with zero seminorm (\cite[Sect.1]{Hand81}).
We note that, in general, $J^*(\mathcal A) \subseteq \Jmin(\mathcal{A}_b)$. 
But if $\ring C = \mathbb{C}$, then in fact $J^*(\mathcal A) = 
\Jmin(\mathcal{A}_b)$, since in this case the $^*$-algebra 
$\mathcal{A}_b/{J^*(\mathcal A)}$ can be completed to a $C^*$-algebra and
hence admits a faithful $^*$-representation.
\end{remark}

%
%

\section{Morita Equivalent $^*$-Algebras}
\label{MoritaSec}

In this section, we shall discuss the behavior of the lattice of
(closed) $^*$-ideals under algebraic Rieffel induction and prove that
the lattice of (closed) $^*$-ideals is a `formal Morita invariant'.

Let $\BXA$ be a $(\mathcal B$-$\mathcal A)$-bimodule
with an $\mathcal A$-valued inner product (see Section \ref{PrelimSec}).
We  define the
map $\PhiX: 2^{\mathcal A} \to 2^{\mathcal B}$ by
\begin{equation}
\label{PhiXDef}
    \PhiX (\mathcal I) := \{B \in \mathcal B \; | \;
                          \forall x, y \in \BXA: 
                          \SPA{x, B\cdot y} \in \mathcal I\}.
\end{equation}
If we consider (closed) $^*$-ideals of $\mathcal A$, then we obtain
the following properties of $\PhiX$:
\begin{proposition}
\label{PhiXAlmostHomoProp}
Let $\mathcal A$, $\mathcal B$ be $^*$-algebras over $\ring C$ and let
$\BXA$ be a $(\mathcal B$-$\mathcal A)$-bimodule with 
$\mathcal A$-valued inner product. Then $\PhiX$ maps $^*$-ideals of 
$\mathcal A$ to $^*$-ideals of $\mathcal B$ and satisfies
\begin{equation}
\label{PhiXAlmostHomo}
    \PhiX(\mathcal I) \wedge \PhiX(\mathcal J) = 
    \PhiX(\mathcal I \wedge \mathcal J)
    \qquad
    \text{and}
    \qquad
    \mathcal I \le \mathcal J \Longrightarrow 
    \PhiX (\mathcal I) \le \PhiX(\mathcal J)
\end{equation}
for all $^*$-ideals $\mathcal I$, $\mathcal J$ of $\mathcal A$.
If in addition the inner product is positive semi-definite and
satisfies \textbf{P}, then $\PhiX$ maps closed $^*$-ideals to closed
$^*$-ideals, satisfies (\ref{PhiXAlmostHomo}) with respect to the
lattice of closed $^*$-ideals, and one has
\begin{equation}
\label{RieffelXPhiX}
    \ker (\RieffelX \piA) = \PhiX (\ker\piA)
\end{equation}
for any $^*$-representation $\piA$ of $\mathcal A$.
\end{proposition}
\begin{proof}
The first part is a simple computation thus assume that the inner
product is positive semi-definite and satisfies \textbf{P}. Let $\piA$
be a $^*$-representation of $\mathcal A$ on $\mathfrak H$ and let
$\piB = \RieffelX\piA$ be the induced $^*$-representation of 
$\mathcal B$. Then $\piB (B) = 0$ if and only if 
$[B \cdot y \otimes \phi] = 0$ for all $y \in \BXA$ and 
$\phi \in \mathfrak H$ since the equivalence classes of those 
$y \otimes \phi$ span $\mathfrak K$. But since the inner product of
$\mathfrak K$ is non-degenerate this is equivalent to 
$0 = \SPK{[x\otimes \psi], [B \cdot y \otimes \phi]} = 
\SPH{\psi, \piA(\SPA{x, B\cdot y})\phi}$ for all $x, y \in \BXA$ and
$\phi, \psi \in \mathfrak H$. Thus it is equivalent to 
$\piA(\SPA{x, B\cdot y}) = 0$ and hence $B \in \PhiX(\ker\piA)$ which
implies (\ref{RieffelXPhiX}). Hence closed $^*$-ideals are mapped to
closed $^*$-ideals and since $\wedge$ and $\le$ are given by
$\cap$ and $\subseteq$ the proof is completed.
\end{proof}

Note that from (\ref{PhiXAlmostHomo}) one obtains that $\PhiX$ is 
almost a homomorphism of lattices, i.e. we have
\begin{equation}
\label{PhiXAlmostHomoIII}
   \PhiX(\mathcal I) \vee \PhiX (\mathcal J)
   \le
   \PhiX (\mathcal I \vee \mathcal J)
\end{equation}
for all $^*$-ideals or closed $^*$-ideals, respectively. Note that in
general we only can guarantee `$\le$' since there may exist more and
perhaps smaller (closed) $^*$-ideals of $\mathcal B$ which are not in
the image of $\PhiX$.

In order to discuss Morita equivalence, we need 
the following technical lemma, which does not use any positivity
requirements but only the $^*$-structures.
\begin{lemma}
\label{TechnicLem}
Let $\mathcal A$, $\mathcal B$ be $^*$-algebras over $\ring C$ and
$\BXA$ a $(\mathcal B$-$\mathcal A)$-bimodule with compatible 
$\mathcal A$- and $\mathcal B$-valued inner products 
and let $\mathcal I \subseteq \mathcal A$ be a
$^*$-ideal. Then $A \in \PhiccX(\PhiX (\mathcal I))$ if and only if
for all $x, y, z, w \in \BXA$ one has
\begin{equation}
\label{PhiXPhiXI}
    \SPA{x, y} A \SPA{z, w} \in \mathcal I.
\end{equation}
\end{lemma}
\begin{proof}
By definition $A \in \PhiccX(\PhiX(\mathcal I))$ if and only if one
has $\SPA{y, \SPB{\cc x, A\cdot\cc z}\cdot w} \in \mathcal I$ for
all $\cc x, \cc z \in \AXB$ and $y, w \in \BXA$. Then it is
a simple computation to obtain the condition (\ref{PhiXPhiXI}).
\end{proof}

If in addition the $\mathcal A$-valued inner product is full, then
condition (\ref{PhiXPhiXI}) means that the ideal generated by $A$ must
be contained in $\mathcal I$. If furthermore $\mathcal A$ has an
approximate identity, then we have 
$\PhiccX(\PhiX(\mathcal I)) = \mathcal I$ for all $^*$-ideals
$\mathcal I$. By symmetry in $\mathcal A$ and $\mathcal B$, we end up
with the following proposition.
\begin{proposition}
\label{IdealLatticeProp}
Let $\mathcal A$, $\mathcal B$ be $^*$-algebras over $\ring C$ with approximate 
identities
and let $\BXA$ be a  $(\mathcal B$-$\mathcal A)$-bimodule with compatible full 
$\mathcal B$- and $\mathcal A$-valued inner products. 
Then the lattices of $^*$-ideals are isomorphic via $\PhiX$
and $\PhiccX$. 
\end{proposition}
\begin{proof}
The fact that $\PhiX$ and $\PhiccX$ are inverse to each other follows
from the above considerations. Hence it remains to show that $\PhiX$
(and thus $\PhiccX$) is a lattice homomorphism. But this follows from
bijectivity and Proposition~\ref{PhiXAlmostHomoProp}.
\end{proof}

In order to make statements about closed $^*$-ideals, we need the
positivity requirements for the inner products to guarantee
that $\RieffelX \piA$ is indeed a $^*$-representation.
\begin{theorem}
\label{ClosedIdealLatticeTheo}
Let $\mathcal A$, $\mathcal B$ be idempotent $^*$-algebras over $\ring C$ 
and let $\BXA$ be an equivalence bimodule. Then
the lattices of closed $^*$-ideals are isomorphic via $\PhiX$ and
$\PhiccX$.  
\end{theorem}
\begin{proof}
Let $\mathcal I = \ker \piA$. We
may assume, due to Lemma~\ref{SNDLem}, that $\piA$ is strongly
non-degenerate. Hence, by \cite[Lem.~5.7]{BuWa99a}, we know that
$\RieffelXc\RieffelX\piA$ is unitarily equivalent to $\piA$. So they
have the same kernel and by (\ref{RieffelXPhiX}) it follows that
$\PhiccX(\PhiX(\mathcal I)) = \mathcal I$. Analogously, one obtains
$\PhiX(\PhiccX(\mathcal J)) = \mathcal J$ for all closed $^*$-ideals
of $\mathcal B$ and hence $\PhiX$ and $\PhiccX$ are inverse to each
other. The homomorphism properties follow as in
Proposition~\ref{IdealLatticeProp}.
\end{proof}

Note that the so-called Rieffel correspondence theorem
(\cite[Thm.~3.24]{RaeWil98}) can be recovered as a consequence of the above
theorem. Indeed, if $\mathcal A$ and $\mathcal B$ are strongly Morita
equivalent $C^*$-algebras, then by Theorem
\ref{TheNiceTheoremAboutCStar} it follows that 
$\mathcal{P}_{\mathcal A}$ and $\mathcal{P}_{\mathcal B}$ are formally
Morita equivalent. So $\mathfrak{L}(\mathcal{P}_{\mathcal A})$ and 
$\mathfrak{L}(\mathcal{P}_{\mathcal B})$ are isomorphic by 
Theorem \ref{ClosedIdealLatticeTheo} and therefore 
$\mathfrak{L}(\mathcal A)$ and $\mathfrak{L}(\mathcal B)$ are
isomorphic by Proposition \ref{SameLatticeProp}. We also observe the following
corollary.
\begin{corollary}\label{TrivialJminCor}
Let $\mathcal A$ and $\mathcal B$ be non-degenerate and idempotent
$^*$-algebras over $\ring C$ which are formally Morita equivalent. 
Then  $\Jmin (\mathcal A)=0$ if and only if 
$\Jmin (\mathcal B)=0$.
\end{corollary}
\begin{proof}
Let $\BXA$ be an equivalence bimodule. Note that if $\mathcal A$ 
is non-degenerate and 
idempotent, then $\BSP{x,y\cdot A}=0$ for all $x,y \in \mathfrak X$ implies
that $A=0$ (see the proof of Lemma \ref{FormalImpliesStarLem}). 
Hence the conclusion follows.
\end{proof}

Let us finally discuss the general phenomena of
$^*$-algebras having equivalent (even isomorphic) categories of
(strongly non-degenerate) $^*$-representations but which are not
formally Morita equivalent. The following discussion generalizes
the example of $\ring C$ and $\bigwedge^\bullet (\ring C^n)$ from
\cite[Cor.~5.20]{BuWa99a}.

Consider the categories of $^*$-representations of
$\mathcal A$ and $\mathcal A \big/ \Jmin(\mathcal A)$. More generally,
we shall consider a $^*$-ideal 
$\mathcal J \subseteq \Jmin (\mathcal A)$ and the $^*$-algebra
$\mathcal A \big/ \mathcal J$. Denoting by 
$\rho: \mathcal A \to \mathcal A \big/ \mathcal J$ the
$^*$-homomorphism $A \mapsto [A]$ we can pull-back
$^*$-representations 
$\pi: \mathcal A \big/ \mathcal J \to \Bounded(\mathfrak H)$ 
by
\begin{equation}
\label{PullBackRep}
    \rho^* \pi (A) := \pi(\rho(A)) = \pi([A]),
\end{equation}
and obtain a $^*$-representation $\rho^*\pi$ of $\mathcal A$ on the
same pre-Hilbert space $\mathfrak H$. This is clearly functorial,
i.e. compatible with intertwiners. Moreover, if $\pi$ is also strongly
non-degenerate then $\rho^*\pi$ is strongly non-degenerate. On
the other hand, let $\pi: \mathcal A \to \Bounded (\mathfrak H)$ be a
$^*$-representation of $\mathcal A$ then 
$\mathcal J \subseteq \Jmin (\mathcal A) \subseteq \ker\pi$. Hence we
can also push-forward the $^*$-representation and obtain a
$^*$-representation $\rho_* \pi$ of $\mathcal A \big/ \mathcal J$
on the same representation space
\begin{equation}
\label{PushForwardRep}
    \rho_*\pi ([A]) := \pi(A).
\end{equation}
Again this is functorial and if $\pi$ is strongly non-degenerate then
$\rho_*\pi$ is also strongly non-degenerate. 
Collecting these results, we obtain 
\begin{proposition}
\label{AQuotientJminProp}
Let $\mathcal A$ be a $^*$-algebra over $\ring C$ and 
$\mathcal J \subseteq \Jmin (\mathcal A)$ a $^*$-ideal. Then the
functors $\rho^*$ and $\rho_*$ implement an isomorphism of the
categories $\srepA$ and $\srepAJ$ as well as of $\sRepA$ and $\sRepAJ$. 
\end{proposition}
We also note that $\mathcal A$ and $\mathcal A \big/ \mathcal J$ have isomorphic
lattices of closed $^*$-ideals, for any $^*$-ideal 
$\mathcal J \subseteq \Jmin (\mathcal A)$. Let
$\mathcal I=\ker \pi \in \mathfrak{L}(\mathcal A)$. 
Since $\mathcal J \subseteq \Jmin (\mathcal A)\subseteq \ker \pi$, it follows
that  
$\rho(\mathcal I)= \ker \rho_*\pi$. Thus $\rho$ induces a map
$\mathfrak{L}(\mathcal A) \to \mathfrak{L}(\mathcal A \big/ \mathcal J)$ which
can be easily checked to be a bijection which preserves ordering, and hence
is a lattice isomorphism. We can then state
\begin{proposition}\label{IsomLatticeProp}
Let $\mathcal A$ be a $^*$-algebra over $\ring C$ and 
$\mathcal J \subseteq \Jmin (\mathcal A)$ a $^*$-ideal. Then the natural
projection $\rho: \mathcal A \to \mathcal A \big/ \mathcal J$ induces an
isomorphism of $\mathfrak{L}(\mathcal A)$ and 
$\mathfrak{L}(\mathcal A \big/ \mathcal J)$. In particular, 
$\Jmin(\mathcal A \big/ \mathcal J)=\Jmin(\mathcal A) \big/ \mathcal J$.
\end{proposition}

We finally observe that this discussion provides a large class of examples of
$^*$-algebras over $\ring C$ having equivalent categories of
$^*$-representations, isomorphic lattice of $^*$-ideals but which are not
formally Morita equivalent.
\begin{proposition}
\label{NotMoritaEquivProp}
Let $\mathcal A$ be a $^*$-algebra over $\ring C$ with approximate
identity and let $\mathcal I$ be a $^*$-ideal contained in
$\Jmin(\mathcal A)$ but $\mathcal I \ne \Jmin (\mathcal A)$. Then
$\mathcal A \big/ \mathcal I$ and $\mathcal A \big/ \Jmin(\mathcal A)$
have isomorphic categories of (strongly non-degenerate)
$^*$-representations and isomorphic lattices of closed $^*$-ideals,
but they are not formally Morita equivalent.
\end{proposition}
\begin{proof}
The first statement follows directly from
Propositions~\ref{AQuotientJminProp} and \ref{IsomLatticeProp} by transitivity. 
It also follows from Proposition \ref{SameLatticeProp} that 
$\mathcal A \big/ \Jmin (\mathcal A)$ has trivial minimal ideal, while
$\Jmin(\mathcal A \big/ \mathcal I)=\Jmin(\mathcal A)/\mathcal I \neq 0$.
So the result follows by Corollary \ref{TrivialJminCor}.
\end{proof}

%
%

\section{Examples from Differential Geometry}
\label{FunctionSec}

In this section, we describe closed $^*$-ideals and minimal ideals
of $^*$-algebras arising in differential geometry.
We start with a description of the closed $^*$-ideals
in $C^\infty(M)$, the algebra of complex-valued smooth functions on an
$m$-dimensional smooth manifold $M$.

Let us first recall that  given a positive
Borel measure $\mu$ on $M$, with compact support, we can define a 
positive linear functional $\omega_\mu$ on $C^\infty(M)$ by integration 
with respect to this measure and consider 
the corresponding GNS representation of
$C^\infty(M)$, denoted by $\pi_\mu$. 
We observe that in this case, $\ker \pi_\mu = 
\{ f \in C^\infty(M) \mbox{ with }\int f |g|^2 d\mu = 0 \mbox{ for all }
g \in C^\infty(M) \}$ coincides with the Gel'fand ideal 
$\{f \in C^\infty(M) \mbox{ with }\int |f|^2 d\mu = 0 \}$.
To each closed set $F \subseteq M$, there naturally corresponds a $^*$-ideal
of $C^\infty(M)$ given by
\begin{equation}
\label{IdealF}
\mathcal{I}_F=\{f \in C^\infty(M) \, 
\mbox{ such that } \, f|_F = 0 \},
\end{equation}
called the \emph{vanishing ideal} of $F$. The following proposition
shows that these ideals are exactly the closed ideals of $C^\infty(M)$.
\begin{proposition}
A $^*$-ideal $\mathcal I \subseteq C^\infty(M)$ is closed if and only if
it is the vanishing ideal of some closed set $F \subseteq M$.
Moreover, $\mathcal I = \ker \pi_{\omega}$, for some GNS representation 
$\pi_\omega$, if and only if $\mathcal{I} = \mathcal{I}_K$ for some
compact set $K \subseteq M$.
\end{proposition}
\begin{proof}
Let $F \ne \emptyset$ be a closed set in $M$. 
Denote by $\mathcal{M}_F$ the set of
positive compactly supported Borel measures on $M$ with support in
$F$. By the continuity of functions in $C^\infty(M)$, it follows that 
$\mathcal I_F = \bigcap_{\mu \in \mathcal{M}_F} \ker \pi_{\mu}$.
Thus $\mathcal
I_F$ is closed. Conversely, if $\mathcal I$ is closed, then it is
the intersection of kernels of GNS representations of some positive
linear functionals. But the positive linear functionals in 
$C^\infty (M)$ are given by positive Borel measures, see
e.g. \cite[App.~B]{BuWa99a}. So we can write $\mathcal I = \ker \pi$, for
$\pi = \oplus_\alpha \pi_{\mu_\alpha}$ and define $F = 
\overline{\bigcup_\alpha \supp \mu_\alpha}$. It is then easy to check that
$\mathcal I = \mathcal{I}_F$.
The last assertion follows from the fact that if $K$ is compact, then
there is a positive Borel measure with support exactly $K$.
\end{proof}

We recall that the vanishing ideals of closed sets in $C^\infty(M)$
are exactly the closed ideals with respect to the (weak)
$C^0$-topology induced from $C(M)$, that is, the locally convex
topology generated by the family of semi-norms 
$\rho_K(f)=\sup_{x \in K}|f(x)| $, for $K \subseteq M$ compact (this
follows from \cite[Thm.~3]{Whit48}). Hence we have 
\begin{corollary}
A $^*$-ideal $\mathcal I$ in $C^\infty(M)$ is closed (in the algebraic
sense) if and only if it is closed with respect to the (weak)
$C^0$-topology.  
\end{corollary}
We remark that similar results hold for the 
algebra $C^\infty_0(M)$ of compactly supported smooth complex-valued
functions on $M$. 
\begin{proposition}
Let $\mathcal I \subseteq C^\infty_0(M)$ be a $^*$-ideal.
Then $\mathcal I$ is closed if and only if $\mathcal I$ is a vanishing
ideal if and only if $\mathcal I$ is the kernel of a GNS representation.
\end{proposition}

It is clear that the minimal ideal of $C^\infty(M)$ is trivial. We
will now consider an example with nontrivial minimal ideal, namely
the sections of the Grassmann algebra bundle over $M$. 
This generalizes the example of the Grassmann algebra
discussed in the introduction. We denote by $\bigwedge^\bullet T^*M$
the complexified Grassmann algebra bundle over $M$. There are
several possibilities to define fibrewise $^*$-involutions for
$\bigwedge^\bullet T^*M$. We do it by declaring the 
real one-forms to be Hermitian
and extending this to an antilinear anti-automorphism of the Grassmann
algebra bundle. Then $\Gamma^\infty (\bigwedge^\bullet T^*M)$ becomes
a $^*$-algebra over $\mathbb C$. Denote by $\bigwedge^+ T^*M$ the
bundle of forms of positive degree. We have the following
characterization of the minimal ideal of 
$\Gamma^\infty (\bigwedge^\bullet T^*M)$.
\begin{proposition}
\label{GrassmannBundleProp}
The minimal ideal of $\Gamma^\infty (\bigwedge^\bullet T^*M)$ is given
by $\Gamma^\infty (\bigwedge^+ T^*M)$. Thus the representation
theories of $\Gamma^\infty (\bigwedge^\bullet T^*M)$ and 
$C^\infty (M)$ are equivalent but the algebras are not formally Morita
equivalent. The closed ideals of 
$\Gamma^\infty (\bigwedge^\bullet T^*M)$ are of the form 
$\mathcal I \oplus \Gamma^\infty (\bigwedge^+ T^*M)$, where 
$\mathcal I$ is a closed ideal of $C^\infty (M)$.
\end{proposition}
\begin{proof}
Due to the nilpotence properties of forms with positive degree, it is
easy to see that $\Gamma^\infty (\bigwedge^+ T^*M)$ has to be
contained in the minimal ideal. On the other hand,
$C^\infty (M) = \Gamma^\infty (\bigwedge^0 T^*M)$ has sufficiently many
positive linear functionals and thus the first statement follows. The
other statement is clear from the above results and
Section~\ref{MoritaSec}. 
\end{proof}

Let $E \rightarrow M$ be a 
complex Hermitian vector bundle over $M$, with $k$-dimensional fibers and
fiber metric $h$. We write $\End(E) \rightarrow M$ for the 
corresponding endomorphism bundle and $\Gamma^\infty(E)$, $\Gamma^\infty(\End(E))$
for the spaces of smooth sections of $E$ and $\End(E)$, respectively.
We denote the trivial bundle $M\times \mathbb{C}^k
\rightarrow M$ by $t(\mathbb{C}^k)$. The final part of this section shows that
$\Gamma^\infty(E)$ is a 
($\Gamma^\infty(\End(E))$-$C^\infty(M)$)-equivalence bimodule and gives
a description of the closed $^*$-ideals in $\Gamma^\infty(\End(E))$.

We recall that $E$ can always be regarded as a subbundle of some
trivial bundle $t(\mathbb{C}^N)$ over $M$ (the proof in 
\cite[Ch.~I, Thm.~6.5]{Karoubi78} works for arbitrary manifolds since in
this case vector bundles are of finite type, see \cite[Sect.~5]{MilSta74}, 
\cite[Lem.~2.7]{Munkres63}). 
Denote $C^\infty(M)$ by $\mathcal A$.
The natural Hermitian metric in $t(\mathbb{C}^N)$ induces an $\mathcal A$-valued
inner product on $\Gamma^\infty(t(\mathbb{C}^N)) \cong \mathcal{A}^N$ and the
embedding $i : E \hookrightarrow t(\mathbb{C}^N)$ can be assumed to preserve
metrics (\cite[Ch.~I, Thm.~8.8]{Karoubi78}). 
Identifying $\Gamma^\infty(E)$ with a submodule of 
$\mathcal{A}^N$, we can write 
$\Gamma^\infty(E) \oplus {\Gamma^\infty(E)}^\perp \cong \mathcal{A}^N$ and
hence there is a projection $Q \in M_N(\mathcal A)$ so that 
$\Gamma^\infty(E)\cong Q\mathcal{A}^N$ as right $\mathcal A$-modules.
Note also that $Q M_N(\mathcal A)Q \cong \Gamma^\infty(\End(E))$ 
as $^*$-algebras. 
It is clear that $\Gamma^\infty(E)$ is a 
($\Gamma^\infty(\End(E))$-$C^\infty(M)$)-bimodule equipped with
$\Gamma^\infty(\End(E))$ and $C^\infty(M)$ valued inner products
given by $(u,v) \mapsto \Theta_{u,v} \in \Gamma^\infty(\End(E))$, where 
$\Theta_{u,v}(w)=u\cdot h(v,w)$, and
$(u,v) \mapsto h(u,v) \in C^\infty(M)$, for $u,v,w \in \Gamma^\infty(E)$.
We finally observe that this bimodule structure on $\Gamma^\infty(E)$
(with the algebra valued inner products) is just the one considered for
$Q\mathcal{A}^N$ in \cite[Sect.~6]{BuWa99a} after the 
aforementioned identifications.
Recall that a projection $Q \in M_N(\mathcal A)$ is called 
\emph{full} if $\mathbb C$-span
$\{TQS \; ; \; T,S \in M_N(\mathcal A)\} = M_N(\mathcal A)$. The following 
observation is well-known.
\begin{lemma}\label{FullProp}
Let $Q \in M_N(C^\infty (M))$ be a non-zero projection. Then $Q$ is full.
\end{lemma}

It then follows from \cite[Sect.~6]{BuWa99a} that $\Gamma^\infty(E)$
satisfies all the properties of a 
($\Gamma^\infty(\End(E))$-$C^\infty(M)$)-equivalence bimodule, except possibly
for property \axiom{Q}, which we now show holds in general. 
Let $(\pi, \mathfrak{H})$ be a $^*$-representation of 
$QM_N(\mathcal A)Q \cong \Gamma^\infty(\End(E))$ and 
$\tildeK = \overline{Q\mathcal{A}^N}\otimes_{QM_N(\mathcal A)Q}\mathfrak{H}$.
Note that $\overline{Q\mathcal{A}^N}$ is again just $\Gamma^\infty(E)$ (with
left and right actions given by the adjoint of the previous ones).
Given $s_i \in \Gamma^\infty(E)$ and $\psi_i \in \mathfrak{H}$, we must show 
that
\begin{equation} \label{positivity}
\SPKT{\sum_i s_i\otimes\psi_i,\sum_j s_j\otimes\psi_j} \geq 0,
\end{equation}
where $\SPKT{s_1\otimes\psi_1,s_2\otimes\psi_2} = 
\SPH{\psi_1, \pi(\SP{s_1,s_2})\psi_2}$.
We then have
\begin{lemma}
$\Gamma^\infty(E)$ satisfies property \axiom{Q}.
\end{lemma}
\begin{proof}
Assuming $E = t(\mathbb{C}^k)$, the $s_i$'s are just 
$\mathbb{C}^k$-valued functions on $M$ and $\End(E) = t(M_k(\mathbb{C}))$.
Note that if $e_j$ denote the canonical section basis of $E$, then for each
$i$ one can find a $T_i \in t(M_k(\mathbb{C}))$ so that $e_1\cdot T_i = s_i$
(recall that $e_1 \cdot T_i$ means ${T_i}^*e_1$, 
where the last expressions is just
the usual matrix multiplication). For instance, define $T_i$ to be the
matrix valued function with $s_i$ as the first row and zero everywhere else.
Now the proof follows the usual trick:
$ \SPKT{\sum_i s_i\otimes \psi_i,\sum_j s_j\otimes \psi_j} =
\SPKT{\sum_i e_1 T_i\otimes \psi_i,\sum_j e_1T_j\otimes \psi_j} =
\SPKT{ e_1 \otimes \sum_i T_i\psi_i, e_1\otimes \sum_j T_j\psi_j} \geq 0
$.

To prove the result in general, let $U_k$, $k=1, \ldots, r$ 
be an open cover of $M$ so that
$E|_{U_k}$ is trivial. Let $\chi_k$ be a quadratic partition 
of unity subordinated to this cover. Then observe that
for $s_1, s_2 \in \Gamma^\infty(E)$, we have that
$\langle s_1, s_2 \rangle_{QM_N(\mathcal A)Q} = \sum_k 
\langle \chi_k s_1, \chi_k s_2 \rangle_{QM_N(\mathcal A)Q}$
(since $\sum_k \SP{\chi_k s_1, \chi_k s_2}(x) = \sum_k 
\Theta_{\chi_k s_1(x), \chi_k s_2(x)} = 
\sum_k \chi_k(x) \overline{\chi_k}(x) \Theta_{s_1(x), s_2(x)} = 
\SP{s_1,s_2}(x)$).
Hence, 
$$\SPKT{\sum_i s_i\otimes \psi_i,\sum_i s_i\otimes \psi_i} = 
\sum_{i,j}\SPKT{\psi_i,\SP{s_i,s_j}\psi_j} =
\sum_{i,j,k} \SPKT{\psi_i,\SP{\chi_ks_i,\chi_ks_j}\psi_j}.
$$
Now note that the last expression  is just
$\sum_k \SP{\sum_i (\chi_ks_i)\otimes\psi_i,\sum_i (\chi_ks_i)\otimes\psi_i}$.
But since $\supp \chi_ks_i \subseteq U_k$, and $E|_{U_k}$ is trivial, the
result follows from the previous discussion.
\end{proof}
\begin{proposition}\label{EquivProp}
If $E$ is a complex hermitian vector bundle over $M$, then $\Gamma^\infty(E)$
is a $(C^\infty(M)$-$\Gamma^\infty(\End(E)))$-equivalence bimodule.
\end{proposition}

One can also consider compactly supported functions and sections. In this case,
$C_0^\infty(M)$ and $\Gamma^\infty_0(\End(E))$ are no longer unital in general,
but have approximate identities. So the results of \cite[Sect.~6]{BuWa99a}
can still be applied and the same proof used for Prop.~\ref{EquivProp} shows that
\begin{proposition}
$\Gamma_0^\infty(E)$
is a $(C_0^\infty(M)$-$\Gamma_0^\infty(\End(E)))$-equivalence bimodule.
\end{proposition}

We finally remark that combining the results of Section \ref{MoritaSec},
one easily arrives at a characterization of the closed $^*$-ideals of
$\Gamma^\infty(\End(E))$ (and analogously for $\Gamma_0^\infty(\End(E))$).
\begin{corollary}
A $^*$-ideal $\mathcal I \subseteq \Gamma^\infty(\End(E))$ 
is closed if and only it is a vanishing ideal, i.e. of the form
$\{B \in \Gamma^\infty(\End(E)) \;|\; B(x) = 0 \; \, \forall x \in F\}$ 
for some
closed $F \subseteq M$.
\end{corollary}

%
%

\section*{Acknowledgements}

We would like to thank Prof.~Ara for his clarifying remarks and for
pointing out to us the paper \cite{Hand81} and Martin Bordemann 
for many useful discussions, in particular
concerning the minimal ideal. We would also like to thank 
Pierre Bieliavsky, Michel Cahen,
Laura DeMarco, Diogo Gomes, Simone Gutt, Daniel Markiewicz, and Alan
Weinstein for helpful remarks and discussions.

%
%

\begin{small}

\end{small}
\end{document}